\documentclass[12pt]{article}
\usepackage{amssymb,amsmath,amsfonts}
\usepackage{stmaryrd}
\usepackage{color}
\usepackage[latin1]{inputenc}
\usepackage{enumitem}

\usepackage{hyperref}

\usepackage[small,bf,hang]{caption} 

\usepackage{epsfig}
\usepackage[font=small,skip=4pt]{caption}

\def\BBox{\kern  -0.2cm\hbox{\vrule width 0.2cm height 0.2cm}}

\newtheorem{theorem}{Theorem}[section]
\newtheorem{lemma}[theorem]{Lemma}

\newtheorem{definition}[theorem]{Definition}
\newtheorem{corollary}[theorem]{Corollary}
\newtheorem{proposition}[theorem]{Proposition}
\newtheorem{remark}[theorem]{Remark}
\newtheorem{observation}[theorem]{Observation}
\newtheorem{conjecture}[theorem]{Conjecture}
\newtheorem{problem}[theorem]{Problem}

\graphicspath{{figures/}}

\newenvironment{dedication}
{\begin{quotation}\begin{center}\begin{em}}
{\end{em}\end{center}\end{quotation}}

\setlist{noitemsep}

\title{\textbf{Colourings of cubic graphs inducing isomorphic monochromatic subgraphs }}

\author{Mari\'{e}n Abreu$^{1}$\thanks{The research that led to the present paper was partially supported by a grant of the group GNSAGA of INdAM and by the Italian Ministry Research Project PRIN 2012 ``Geometric Structures, Combinatorics and their Applications''.}, Jan Goedgebeur$^2$\thanks{Supported by a Postdoctoral Fellowship of the Research Foundation Flanders (FWO).
\newline {\em Email addresses:} marien.abreu@unibas.it (M. Abreu),~ jan.goedgebeur@ugent.be (J. Goedgebeur), ~ domenico.labbate@unibas.it (D. Labbate), ~ giuseppe.mazzuoccolo@univr.it (G. Mazzuoccolo)}, \\ Domenico Labbate$^{1}$$^*$, Giuseppe Mazzuoccolo$^{3}$$^*$
\\[2ex]
\footnotesize $^1$ Dipartimento di Matematica, Universit\`{a} degli Studi della Basilicata,\\[-3mm]
\footnotesize Viale dell'Ateneo Lucano, I-85100 Potenza, Italy.\\
\footnotesize $^2$ Deparment of Applied Mathematics, Computer Science and Statistics\\[-3mm]
\footnotesize Ghent University, Krijgslaan 281 - S9, 9000 Ghent, Belgium.\\
\footnotesize $^3$ Dipartimento di Informatica,\\[-3mm]
\footnotesize Universit\`{a} degli Studi di Verona, Strada le Grazie 15, 37134 Verona, Italy.}

\date{}

\begin{document}
\maketitle

\begin{dedication}
In loving memory of Ella.
\end{dedication}

\begin{abstract}
A \emph{$k$--bisection} of a bridgeless cubic graph $G$ is a $2$--colouring of its vertex set such that the colour classes have the same cardinality and all connected components in the two subgraphs induced by the colour classes ({\em monochromatic components} in what follows) have order at most $k$.
Ban and Linial conjectured that {\em every bridgeless cubic graph admits a $2$--bisection except for the Petersen graph}.
A similar problem for the edge set of cubic graphs has been studied: Wormald conjectured that {\em every cubic graph $G$ with $|E(G)| \equiv 0  \pmod  2$ has a $2$--edge colouring such that the two monochromatic subgraphs are isomorphic linear forests} (i.e.\ a forest whose components are paths).
Finally, Ando conjectured that {\em every cubic graph admits a bisection such that the two induced monochromatic subgraphs are isomorphic.}

In this paper, we provide evidence of a strong relation of the conjectures of Ban--Linial and Wormald with Ando's conjecture. Furthermore, we also give computational and theoretical evidence in their support. As a result, we pose some open problems stronger than the above mentioned conjectures. Moreover, we prove Ban--Linial's conjecture for cubic cycle permutation graphs.

As a by--product of studying $2$--edge colourings of cubic graphs having linear forests as monochromatic components, we also give a negative answer to a problem posed by Jackson and Wormald about certain decompositions of cubic graphs into linear forests.
\end{abstract}



\section{Introduction}\label{intro}

All graphs considered in this paper are finite and simple (without loops or multiple edges). Most of our terminology is standard; for further definitions and notation not explicitly stated in the paper, please refer to~\cite{BM}.

A {\em bisection} of a cubic graph $G=(V,E)$ is a partition of its vertex set $V$ into two disjoint subsets $({\cal B}, {\cal W})$ of the same cardinality. Throughout the rest of the paper, we will identify a bisection of $G$ with the vertex colouring of $G$ with two colours $B$ (black) and $W$ (white), such that every vertex of ${\cal B}$ and ${\cal W}$ has colour $B$ and $W$, respectively.  Note that the colourings not necessarily have to be proper and in what follows we refer to a connected component of the subgraphs induced by a colour class as a \emph{monochromatic component}.
Following the terminology introduced in~\cite{BL16}, we define a \emph{$k$--bisection} of a graph $G$ as a $2$--colouring $c$ of the vertex set $V(G)$ such that:

\begin{enumerate}[label=(\arabic*)]
\item[(i)]
$|{\cal B}|=|{\cal W}|$ (i.e.\ it is a bisection), and
\item[(ii)]
each monochromatic component has at most $k$ vertices.
\end{enumerate}

Note that this definition is almost equivalent to the notion of {\em $(k+2)$--weak bisections} introduced by Esperet, Tarsi and the fourth author in~\cite{EMT15}, where there is the further assumption that each monochromatic component is a tree. In particular, note that for $k=2$ the definitions of $2$-bisection and $4$-weak bisection are exactly equivalent. The main reason of this gap in the terminology is due to the fact that in~\cite{BL16}, the authors focus their attention on the order of the monochromatic component while in~\cite{EMT15} the relation with nowhere--zero flows has more focus. Indeed, the existence of a $(k+2)$--weak bisection is a necessary condition for the existence of a nowhere--zero $(k+2)$--flow, as it directly follows from a result by Jaeger~\cite{Jaeger} stating that {\em every cubic graph $G$ with a circular nowhere--zero $r$--flow has a $\lfloor r \rfloor$--weak bisection} (cf.~\cite{Jaeger} and~\cite{EMT15} for more details).

There are several papers in literature considering $2$--colourings of regular graphs which satisfy condition $(ii)$, but not necessarily condition $(i)$, see~\cite{ADOV, BS, HST, LMST}. In particular, it is easy to see that every cubic graph has a $2$--colouring where all monochromatic connected components are of order at most $2$, but, in general, such a colouring does not satisfy condition $(i)$ and so it is not a $2$--bisection. Thus, the existence of a $2$--bisection in a cubic graph is not guaranteed. For instance, the Petersen graph does not admit a $2$--bisection. However, the Petersen graph is an exception since it is the only known bridgeless cubic graph without a $2$--bisection. This led Ban and Linial to state the following conjecture:

\begin{conjecture}[Ban--Linial~\cite{BL16}]\label{banlinial}
Every bridgeless cubic graph admits a $2$--bisection, except for the Petersen graph.
\end{conjecture}

A similar problem for the edge set of a cubic graph has been studied. A {\em linear forest} is a forest whose components are paths and a {\em linear partition} of a graph $G$ is a partition of its edge set into linear forests (cf.~\cite{FTVW,JakWor}). Wormald made the following conjecture:

\begin{conjecture}[Wormald~\cite{Wor}]\label{Wconj}
Let $G$ be a cubic graph with $|E(G)|$ $\equiv 0 \pmod 2$ (or equivalently $|V(G)|$ $\equiv 0 \pmod 4$). Then there exists a linear partition of $G$ into two isomorphic linear forests.
\end{conjecture}

The main aim of this paper is to provide evidence of a strong relation of the two previous conjectures with the following conjecture proposed by Ando~\cite{A17} in the 1990's.

\begin{conjecture}[Ando~\cite{A17}]\label{ando}
Every cubic graph admits a bisection such that the two induced monochromatic subgraphs are isomorphic.
\end{conjecture}

We will use two different approaches: the first one dealing with $2$--bisections and the second one with linear forests.

Generally speaking, Ban--Linial's Conjecture asks for a $2$--vertex colouring of a cubic graph with certain restrictions on the induced monochromatic components. These restrictions deal with the order of the monochromatic components. Hence they also have consequences on the structure of the admissible components. In particular, it is easy to prove that the two induced monochromatic subgraphs in a $2$-bisection are isomorphic. As far as we know, apart the present paper, only Esperet, Tarsi and the fourth author have dealt with the Ban--Linial Conjecture in~\cite{EMT15}. Their main result (see Theorem 11 in~\cite{EMT15}) implies that every (not necessarily bridgeless) cubic graph except the Petersen graph admits a $3$--bisection (cf.\ Section~\ref{blconj}). Moreover, we have proved this conjecture for claw--free cubic graphs in~\cite{AGLM}.

As far as we know, Wormald's Conjecture is still widely open. For some results about it, see e.g.~\cite{FTVW}. Note that the statement of Wormald's Conjecture is equivalent to the statement that {\em there exists a $2$--edge colouring of $E(G)$ such that the two monochromatic subgraphs are isomorphic linear forests}. (In Section~\ref{strongwormaldconj}, we will define such a $2$--edge colouring to be a {\em Wormald colouring}). Also in this case we are dealing with a colouring problem, an edge colouring problem this time, with some restrictions on the admissible monochromatic components. In analogy with the Ban--Linial Conjecture, we will propose a stronger version of Conjecture~\ref{Wconj} by adding a restriction on the order of the monochromatic components (see Problem~\ref{pro:4linearforests}).

We would also like to point out, as mentioned earlier, that we can give a positive answer for some instances of Ando's Conjecture using either $2$--bisections or linear forests. We find it interesting that the few possible {\em exceptions} for which the approach with $2$--bisections is not feasible, are in the class of {\em not} $3$--edge colourable cubic graphs. On the other hand, all known exceptions for which the approach with our strengthened version of Wormald colourings is not feasible, are in the class of $3$--edge colourable cubic graphs. Note that, this indicates a promising possibility for a complete proof of Ando's Conjecture by proving that a given cubic graph cannot be an exception for both of these two approaches (at least for a bridgeless cubic graph).

%



The paper is organised as follows. In Section~\ref{blconj}, we deal with the Ban--Linial Conjecture (cf.\ Conjecture~\ref{banlinial}). First of all, we observe that Ban--Linial's Conjecture implies Ando's Conjecture (cf.\ Proposition~\ref{BanLinialimpliesAndo}). Moreover, in Section~\ref{cyperm}, we prove Ban--Linial's Conjecture for bridgeless cubic graphs in the special case of cycle permutation graphs. In Section~\ref{1connectedBanLinial} we verify that all known cubic graphs without a 2-bisection are not a counterexample for Ando's Conjecture and in Section~\ref{subsect:strongando} we present a possible stronger version of this conjecture.

In Section~\ref{strongwormaldconj}, we propose a strengthened version of Wormald colourings and we show that Ando's Conjecture is strongly related to it (see Proposition~\ref{swortosandocol}).

In Section~\ref{CompRes} we present our computational results. It is known that there are infinitely many counterexamples to Ban--Linial's Conjecture for 1--connected cubic graphs (see~\cite{EMT16}) and in Section~\ref{comput:bl} we produce the complete list of all cubic graphs up to 32 vertices which do not admit a $2$--bisection. In Section~\ref{comput:ando} we show that Ando's conjecture does not have any counterexamples with less than 34 vertices and in Section~\ref{comput:wormald} we present our computational results on Wormald's conjecture.

Finally, in Section~\ref{thomassenprob}, as a by--product of our study of $2$--edge colourings of cubic graphs having linear forests as monochromatic components,  we obtain a negative answer (cf.\ Theorem~\ref{Th:Heawood}) to a problem about {\em linear arboricity} posed by Jackson and Wormald in~\cite{JakWor} (cf.\ also~\cite{Tho}). The linear arboricity of a graph $G$ is the minimum number of forests whose connected components are paths of length at most $k$ required to partition the edge set of a given graph $G$ (cf.\ Section~\ref{thomassenprob} for details). We conclude the paper by posing new problems and conjectures on this topic.

\section{Ban--Linial's Conjecture}\label{blconj}

We recall that a $2$--bisection is a bisection $(B,W)$ such that the connected components in the two induced subgraphs $G[B]$ and $G[W]$ have order at most two. That is: every induced subgraph is a union of isolated vertices and isolated edges.

In~\cite{BL16}, Ban and Linial posed Conjecture~\ref{banlinial} and proved the following Theorem~\ref{BLthm}. Note that the original formulation is in terms of \emph{external bisections}, which are indeed equivalent to $2$--bisections for cubic graphs.

\begin{theorem}\cite[Proposition 7]{BL16}\label{BLthm}
Every $3$--edge colourable cubic graph admits a $2$--bisection.
\end{theorem}

A stronger version of Theorem~\ref{BLthm} is proved in~\cite{EMT16} by Esperet, Tarsi and the fourth author (in terms of $(k+2)$--weak bisections). In particular, they prove the following:

\begin{theorem}\cite{EMT16}\label{thm:circflownr}
If a bridgeless cubic graph $G$ does not admit a $2$--bisection, then it has circular flow number at least $5$.
\end{theorem}

%

%

Because of Theorem~\ref{thm:circflownr}, a possible counterexample to Conjecture~\ref{banlinial} has to be in the class of bridgeless cubic graphs with circular flow number at least $5$ (exactly $5$ if Tutte's famous conjecture is true). An important subclass of bridgeless cubic graphs which are {\em not} $3$--edge colourable  is the class of {\em snarks} (i.e.\ a subclass of not $3$--edge colourable graphs with additional restrictions on the girth and connectivity). There is a vast literature on snarks and their properties -- see, e.g.,~\cite{ALRS, AKLM, BGHM, KaRa10, Ko96-1, MaRa06}. Snarks are especially interesting as they are often smallest possible counterexamples for conjectures: for many graph theory conjectures it can be proved that if the conjecture is false, the smallest counterexample must be a snark. The interested reader might find more information about snarks in, e.g.,~\cite{HS} or~\cite{Zh97}. Snarks with circular flow number $5$ are known (see~\cite{AKLM,EMT16,GMM,MaRa06}) and no snark without a $2$--bisection, other than the Petersen graph, has been found until now.

\begin{figure}[h]
	\centering
		\includegraphics[height=5cm]{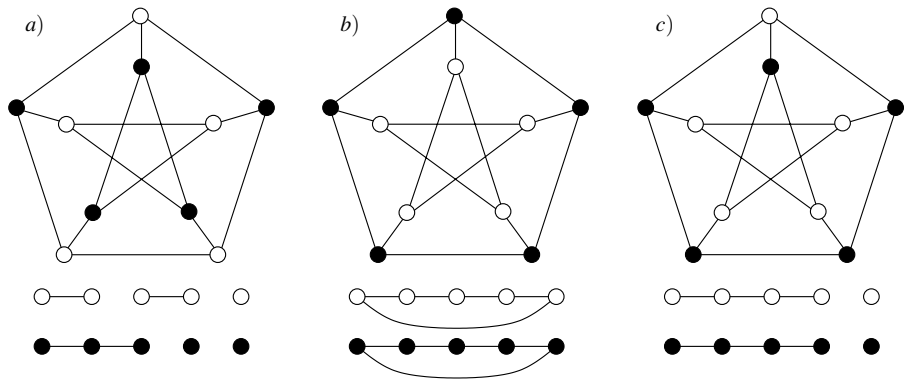}
		\caption{a) A $3$--bisection of $P$; b) and c) Bisections of $P$ with isomorphic parts}\label{figure1} 		
\end{figure}

Conjecture~\ref{banlinial} is strongly related to Ando's Conjecture (i.e.\ Conjecture~\ref{ando}) and the following holds:

\begin{proposition}\label{BanLinialimpliesAndo}
Let $G$ be a bridgeless cubic graph. If $G$ is a counterexample for Ando's Conjecture, then it is also a counterexample for the Ban--Linial Conjecture.
\end{proposition}
\begin{proof}
An easy counting argument shows that the number of isolated edges (vertices) in $G[B]$ and $G[W]$ is the same for every $2$--bisection $(B,W)$ of $G$. Thus every $2$--bisection is a bisection with isomorphic parts. Moreover, it is not hard  to prove that the Petersen graph $P$ is not a counterexample for Ando's Conjecture (see Figure~\ref{figure1}(b) and~\ref{figure1}(c)).
\end{proof}

The previous proposition implies that every cubic graph with a $2$--bisection is not a counterexample for Ando's Conjecture. Not all cubic graphs admit a $2$--bisection: next to the Petersen graph, an infinite family of $1$--connected examples is known~\cite{EMT15}. However, it is proved in~\cite{EMT15} that every (not necessarily bridgeless) cubic graph admits a $3$--bisection. It is not hard  to prove that the Petersen graph has no $3$--bisection with isomorphic parts (see Figure~\ref{figure1}(a)), so there is no chance to prove in general that every cubic graph has a $3$--bisection with isomorphic parts.

In Section~\ref{cyperm} we focus our attention on bridgeless cubic graphs. In particular, we consider a natural generalization of the Petersen graph, namely cycle permutation graphs, and we verify the Ban--Linial Conjecture (and by Proposition~\ref{BanLinialimpliesAndo} also Ando's Conjecture) for all members of this infinite family.
Note that infinitely many of them are not $3$--edge colourable and several authors (see for instance~\cite{BGHM} and~\cite{HagHof}) refer to a member of this subclass as a permutation snark.

In Section~\ref{1connectedBanLinial} we consider the analogous relation between $2$--bisections and bisections with isomorphic parts for the class of $1$--connected cubic graphs. An infinite family of $1$--connected cubic graphs which do not admit a $2$--bisection is known~\cite{EMT15} and it could be that a member of this family might be a counterexample for Ando's Conjecture. In Proposition~\ref{notcounterando} we prove that this is not the case by constructing a bisection with isomorphic parts for all of these graphs. Using a computer program, we also construct all ($1$--connected) cubic graphs of order at most $32$ without a $2$--bisection and we verify that none of them is a counterexample for Ando's Conjecture. We refer the reader to Section~\ref{CompRes} for more details about these computational results.

Finally, in Section~\ref{subsect:strongando} we propose a stronger version of Ando's Conjecture which we will relate to Conjecture~\ref{Wconj} in Section~\ref{strongwormaldconj}.

\subsection{$2$--bisections of cycle permutation graphs}\label{cyperm}

In this section we will introduce cycle permutation graphs and prove that every such graph has a $2$--bisection, except the Petersen graph.

Given a graph $G$ with $n$ vertices labelled $0,1,\ldots,n-1, \, n \ge 4$, and a permutation $\alpha \in S_n$, the symmetric group on the $n$ symbols $\{0,1,\ldots,n-1\}$. The $\alpha$--\emph{permutation graph} of $G$, $P_{\alpha}(G)$ consists of two disjoint copies of $G$, $G_1$ and $G_2$ along with the $n$ edges obtained by joining $i$ in $G_1$ with $\alpha(i)$ in $G_2$, for $i=0,1,\ldots,n-1$. If the graph $G$ is an $n$--cycle labeled $0,1,2,\ldots,n-1$ where $i$ is joined to $i+1$ and $i-1$ $(mod \, n)$ then $P_{\alpha}(G)$ is called a \emph{cycle permutation graph} and will be denoted $C(n,\alpha)$. So, equivalently, a \emph{cycle permutation graph} is a cubic graph consisting of two $n$--cycles $C_1$ and $C_2$, labeled $0, \ldots, n-1$ and such that $i$ in $C_1$ is adjacent to $p(i)$ in $C_2$, for some permutation $p \in S_n$ (where $S_n$ denotes the symmetric group on $\{0,1,\ldots,n-1\}$).

Permutation graphs were first introduced by Chartrand and Harary in 1967~\cite{CH67}, and cycle permutation graphs were given this name in~\cite{Ri79}, but can also be found in~\cite{Ri84} and other references. Clearly, the Petersen graph is a cycle permutation graph and generalised Petersen graphs are cycle permutation graphs as well.

Let $G=C(n,p)$ be a cycle permutation graph, for some permutation $p \in S_n$, having cycles $C_1$ (external) and $C_2$ (internal). We will use the notation $p(i)=p_i$ for all $i \in \{0, \ldots, n-1\}$, as long as subindices are clear, but sometimes we will keep $p(i)$ to avoid triple or quadruple subindices and/or use extra parenthesis to avoid confusion.
Throughout the paper, we assume that $p_0=0$ and we fix the following labelling on the vertices of $G$. Here and along the rest of the paper all indices are taken modulo $n$.

\begin{figure}[h]
	\centering
		\includegraphics[height=5cm]{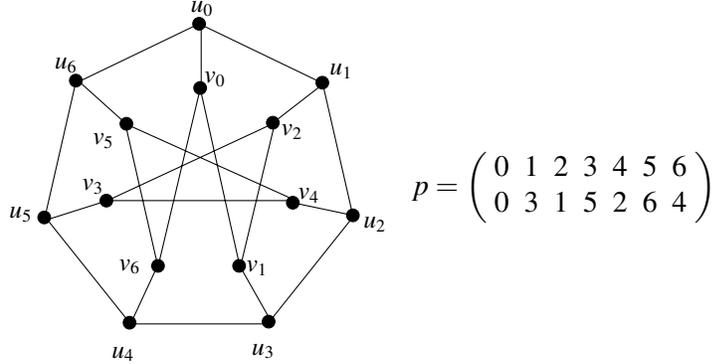}
		\caption{Labelling for a cycle permutation graph}\label{cpglab}
\end{figure}

\begin{itemize}
\item The vertices of the cycle $C_1$ are $u_0, \ldots, u_{n-1}$ with $u_iu_{i+1} \in E(G)$; 
\item The vertices of the cycle $C_2$ are $v_0, \ldots, v_{n-1}$ with $v_iv_{i+1} \in E(G)$; 
\item The edges given by the permutation $p$ are $v_iu_{p_i}$ (see Figure~\ref{cpglab}).
\end{itemize}

Moreover, we always consider the class representatives in $\{0,\ldots,n-1\}$.
For every $p \in S_n$, we denote by $p^{-1}$ the inverse of $p$ and by $\bar{p}$ the permutation such that $\bar{p_i}=n-p_i$  modulo $n$, for all $i$.

\begin{proposition} \label{iso}
The graphs $C(n,p)$, $C(n,p^{-1})$ and $C(n,\bar{p})$ are isomorphic.
\end{proposition}
\begin{proof}
The isomorphism between $C(n,p)$ and  $C(n,p^{-1})$ is easily obtained by switching the roles of $C_1$ and $C_2$.
Moreover, $C(n,\bar{p})$ is obtained from $C(n,p)$ by labeling $C_1$ in counterclockwise order.
\end{proof}

\begin{proposition}\label{neven}
Let $G=C(n,p)$ be a cycle permutation graph. If $n$ is even, then $G$ admits a $2$--bisection.
\end{proposition}

\begin{proof}
For an even $n$, we can alternately colour the edges of $C_1$ and $C_2$ with two colours and the edges not in $C_1$ or $C_2$ with a third colour, then the graph $G$ is $3$--edge colourable, and by Theorem~\ref{BLthm} $G$ admits a $2$--bisection.
\end{proof}

The main purpose of this section is to prove the existence of a $2$--bisection for every cycle permutation graph, except for the Petersen graph (see Theorem~\ref{main}). By Proposition~\ref{neven}, from now on, we can assume that $n$ is odd.
\smallskip

We denote by $||p_{i+1}-p_{i}||$ the length of the path $(u_{p_i},u_{p_i+1},\ldots,u_{p_{i+1}})$ on $C_1$. Note that following such a notation,  $||p_i - p_{i+1}||$ will denote the length of the other path on $C_1$ connecting $p_i$ and $p_{i+1}$, that is: $(u_{p_{i+1}},u_{p_{i+1}+1},\ldots,u_{p_{i}})$.
\smallskip

If $||p_i - p_{i+1}||=1$ or $||p_{i+1}-p_i||=1$  for some $i$, then $G$ is Hamiltonian, and thus $3$--edge colourable and hence admits a $2$--bisection (by Theorem~\ref{BLthm}). Thus, we can also assume that $||p_i - p_{i+1}||> 1$ and $||p_{i+1}-p_i||> 1$ for all $i$.
Finally, since $n$ is odd, $||p_1-p_0||$ and $||\bar{p_1}-\bar{p_0}||$ have different parity, hence we can select $p$ such that $||p_1-p_0||$ is even, that is: $p_1$ is an even integer since $p_0=0$.
We summarise all previous observations in the following remark:

\begin{remark}
To prove the existence of a $2$--bisection for every cycle permutation graph, it is sufficient to prove it for all cycle permutation graphs $C(n,p)$ such that:
$n$ odd,  $p_1$ even, $||p_i - p_{i+1}|| > 1$ and $||p_{i+1}-p_i|| > 1$ for all $i$.
\end{remark}

Let $P$ be a path in a graph $G$, and let $x,y \in V(P)$. We denote by $P(x,y)$ the subpath of $P$ between $x$ and $y$. Let $C$ be a cycle in a graph $G$, and let $x,y,z \in V(C)$. We denote by $C(x,z,y)$ and $C(x,\widehat{z},y)$ the $(x,y)$--path in $C$ containing $z$ and not containing $z$, respectively.

Now, we furnish a sufficient (but in general not necessary) condition to have a $2$--bisection of $G=C(n,p)$. More precisely, we will prove that this condition guarantees the existence of a particular $2$--bisection for $G=C(n,p)$ which satisfies all of the following further properties:

\begin{itemize}
\item The cycle  $C_1$ is alternately black and white, except for two consecutive white vertices, say $a_r$ and $a_{r+1}$;
\item The cycle  $C_2$ is alternately black and white, except for two consecutive black vertices, say $b_s$ and $b_{s+1}$;
\item $a_{r+1}b_s \in E(G)$.
\end{itemize}

Note that the selection of $a_ra_{r+1}$ and $b_sb_{s+1}$ uniquely induces the entire $2$--colouring, which is indeed a bisection.
In order for such a $2$--colouring to have monochromatic components of order at most two, it is enough to find the following {\em good} configuration in~$G$.

\begin{figure}[h]	
	\centering
		\includegraphics[height=4cm]{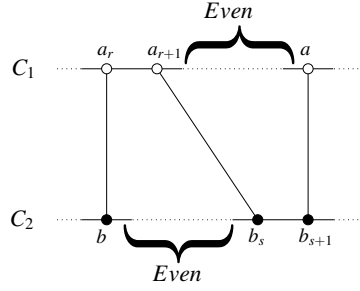}
		\caption{A good configuration}\label{fig:goodconf}		
\end{figure}

\begin{definition}
Let $G=C(n,p)$ be a cycle permutation graph with external cycle $C_1$ and internal cycle $C_2$. A {\em good configuration} in $G$ is a set of six vertices $\{a_r,a_{r+1},a,b_s,b_{s+1},b\} \in V(G)$ such that:
$a_r,a_{r+1},a \in V(C_1)$; $b_s,b_{s+1},b \in V(C_2)$; $a_ra_{r+1}, a_rb, a_{r+1}b_s, ab_{s+1}, b_sb_{s+1} \in E(G)$;
and the paths $C_1(a,\widehat{a_r},a_{r+1})$ (and thus also $C_1(a,\widehat{a_{r+1}},a_r)$) and $C_2(b,\widehat{b_{s+1}},b_s)$ (and thus also $C_2(b,\widehat{b_{s}},b_{s+1})$)  are both of even length (see Figure~\ref{fig:goodconf}).
\end{definition}

\begin{proposition}\label{pro:goodconf}
Let $G=C(n,p)$ be a cycle permutation graph, with $n$ odd. If $G$ has a good configuration, then it admits a $2$--bisection.
\end{proposition}

\begin{proof}
Let  $\{a_r,a_{r+1},a,b_s,b_{s+1},b\}$ be a good configuration of $G$. Set the vertex colour of $a_r,a_{r+1}$ and $a$ white;  the vertex colour of $b_s,b_{s+1}$ and $b$ black; and all other vertices of $G$, alternately black and white in $C_1$ and $C_2$. Note that such alternation is possible and uniquely determined by the previous choices, because the $C_1(a_{r+1},\widehat{a_r},a)$, $C_1(a,\widehat{a_{r+1}},a_r)$, $C_2(b,\widehat{b_{s+1}},b_s)$ and $C_2(b_{s+1},\widehat{b_s},b)$ are all of even length.
All neighbours of the white vertices $a_r,a_{r+1}$ are black, and those of the black vertices $b_s,b_{s+1}$ are white. Moreover, if some edge between vertices of $C_1$ and $C_2$ is monochromatic, by the alternation of colours in those cycles, all of its neighbours will be of the opposite colour. Hence, the monochromatic components have order at most two, that is the $2$--colouring is a $2$--bisection of $C(n,p)$.
\end{proof}

\begin{lemma}\label{notalleven}
Let $G=C(n,p)$ be a cycle permutation graph, with $n$ odd. If there exists an index $i$ such that $||p_{i+1} - p_i|| \equiv 1 \pmod 2$ (lower indices taken modulo $n$), then $G$ admits a $2$--bisection.
\end{lemma}

\begin{proof}
Recall that we are assuming that $p_1$ is even, hence $||p_1 - p_0|| = p_1$ is even. First consider that there is an index $i$ such that $||p_{i+1} - p_i||$ is odd.
Let $t$ be the smallest such index.

In order to find a good configuration in $G$, we consider the path $v_{t-1}v_tv_{t+1}$ in $C_2$ and their neighbours $u_{p_{t-1}}$,$u_{p_{t}}$ and $u_{p_{t+1}}$ in $C_1$, i.e.\ $v_{t-1}u_{p_{t-1}}, v_tu_{p_{t}}, v_{t+1}u_{p_{t+1}} \in E(G)$, $u_{p_{t-1}},u_{p_{t}},u_{p_{t+1}} \in V(C_1)$.

Set $P=C_1(u_{p_{t-1}}, u_{p_{t}}, u_{p_{t+1}})$. We have that $P(u_{p_{t-1}},u_{p_{t}})$ has even length (because $||p_t-p_{t-1}||$ is even) and $P(u_{p_{t}},u_{p_{t+1}})$ has odd length (because $||p_{t+1}-p_t||$ is odd) and greater than $1$ (because $||p_{i+1} - p_i|| > 1$ and $||p_{i} - p_{i+1}|| > 1$  for all $i \in \{0, \ldots, n-1\}$). Let $u$ be the neighbour of $u_{p_{t}}$ in $P(u_{p_{t}},u_{p_{t+1}})$, $v$ be the neighbour of $u$ in $C_2$ and denote by $Q$ the path $C_2(v_{t+1},v_t,v)$. There are two cases:

{\sc Case 1 - LENGTH OF Q ODD:}
Let $b_{s+1}=v_{t-1}$, $b_s=v_t$, $b=v$, $a_{r}=u$, $a_{r+1}=u_{p_t}$ and $a=u_{p_{t-1}}$. Then, $C_2(b,\widehat{b_{s+1}},b_s)$ and $C_1(a,\widehat{a_{r+1}},a_r)$ are both of even length and we have the desired good configuration.

{\sc Case 2 - LENGTH OF Q EVEN:}
Let $b_{s}=v_{t}$, $b_{s+1}=v_{t+1}$, $b=v$, $a_{r+1}=u_{p_{t}}$, $a_{r}=u$ and $a=u_{p_{t+1}}$. Then, $C_2(b,\widehat{b_{s+1}},b_s)$ and $C_1(a,\widehat{a_r},a_{r+1})$ are both of even length and we have the desired good configuration.

In both cases we have a good configuration, so by Proposition~\ref{pro:goodconf} $G$ admits a $2$--bisection.
\end{proof}

Now, we can focus our attention on what happens when we consider the inverse permutation $p^{-1}$.
But first, we collect a list of properties that we are going to use implicitly in the following proofs.

\begin{remark}\label{rempropi&ii}

\begin{enumerate}

\item[(i)] Let $G=C(n,p)$ be a cycle permutation graph  with $n$ odd and $||p_{i+1} - p_{i}||$ even for all $i$. Then, given two consecutive vertices of the inner cycle $v_i, v_{i+1} \in C_2$, we have that for $v_iu_{p_i}, v_{i+1}u_{p_{i+1}} \in E(G)$, the path $(u_{p_i}u_{(p_i)+1}u_{(p_i)+2}\ldots u_{p_{(i + 1)}})$ has even length.

\item[(ii)]  Let $G=C(n,p)$ be a cycle permutation graph  with $n$ odd and $||p^{-1}_{i+1} - p^{-1}_i||$ odd for all $i$. Then, given two consecutive vertices of the outer cycle $u_{p_i}, u_{(p_i)+1} \in C_1$, we have that for $u_{p_i}v_i, u_{p_{i+1}}v_{p^-1((p_i)+1)} \in E(G)$, the path $(v_iv_{i+1}\ldots v_{p^{-1}((p_i)+1)})$ has odd length.

\end{enumerate}
\end{remark}

Recall that if two consecutive vertices of one of the two cycles $C_1$ and $C_2$ have the same colour, we need that their neighbours in the other cycle have the opposite colour.

\begin{lemma}\label{inverseallodd}
Let $G=C(n,p)$, $n$ odd, be a cycle permutation graph. If $||p_{i+1} - p_{i}||$ is even for all $i$ and there exists an index $t$ such that $||p^{-1}_{t+1} - p^{-1}_{t}||$ is even, then $G$ admits a $2$--bisection.
\end{lemma}
\begin{proof}
Again, we make use of Proposition~\ref{pro:goodconf}. Hence, we must only prove the existence of a good configuration.
Consider the vertices $u_{t},u_{t+1}$ in $C_1$ and their neighbours $v_{p^{-1}_{t}}$,$v_{p^{-1}_{t+1}}$ in $C_2$, respectively. Denote by $v$ the vertex adjacent to $v_{p^{-1}_{t}}$ in $C_2$ but not in the $(v_{p^{-1}_{t}},v_{p^{-1}_{t+1}})$--path of $C_2$ of even length. Finally, denote by $u$ the neighbour of $v$ in $C_1$. The $(u,u_t)$--path in $C_2$ of even length does not contain $u_{t+1}$ because of the assumption that $||p_{i+1} - p_{i}||$ is even for all $i$.
Let $b_s=v_{p^{-1}_{t}}$, $b_{s+1}=v$, $b=v_{p^{-1}_{t+1}}$, $a_r=u_t$, $a_{r+1}=u_{t+1}$ and $a=u$. Then, by our assumption on $p$ and $p^{-1}$, $C_2(b,\widehat{b_{s+1}},b_s)$ and $C_1(a,\widehat{a_{r+1}},a_r)$ are both even and we have the desired good configuration.
\end{proof}

In view of the previous lemma, from now on we consider a cycle permutation graph $C(n,p)$ with $n$ odd, $||p_{i+1} - p_{i}||$ even for all $i$ and $||p^{-1}_{i+1} - p^{-1}_{i}||$ odd for all $i$. Note that the Petersen graph lies within this class.

Recall that a Generalised Petersen graph $GP(n,k)$, for $2 \le 2k < n$, 
is a simple cubic graph with vertex
set $V(GP(n,k))=\{u_0, \ldots, u_{n-1}, v_0, \ldots, v_{n-1}\}$ and edge set
$E(GP(n,k))=\{ u_iu_{i+1}, u_iv_i, v_iv_{i+k} : i \in \{0, \ldots, n-1\} \}$ (where all indices are taken modulo $n$).



Even though it is not hard to check that some cycle permutation graphs do not admit a good configuration (e.g.\ the Generalised Petersen graph $GP(9,2)$), we prove in the following theorem that all of them, except for the Petersen graph, admit a $2$--bisection.

\begin{theorem}\label{main}
Let $G=C(n,p)$ be a cycle permutation graph, which is not the Petersen graph. Then $G$ admits a $2$--bisection.
\end{theorem}

\begin{proof} By Proposition~\ref{neven} we can assume that $n$ is odd. By Lemma~\ref{notalleven} we can assume that $p$ satisfies $||p_{i+1} - p_{i}||$ even for all $i$. Finally, by Lemma~\ref{inverseallodd} we can assume that $p^{-1}$ satisfies $||p^{-1}_{i+1} - p^{-1}_{i}||$ odd for all $i$.

First of all, we consider the Generalised Petersen graph $G = GP(n,\frac{n-1}{2})$. If $n=5$ we have the Petersen graph. If $n > 5$, $G$ is $3$--edge colourable (see~\cite{CP}) and by Theorem~\ref{BLthm}, $G$ admits a $2$--bisection.

Hence, we can assume that $p$ satisfies all previous conditions, but it is not $GP(n,\frac{n-1}{2})$.

Let $x,y \in V(G)$ be such that:


We may assume without loss of generality that $i=0$, i.e. $x=v_{p^{-1}(1)}$ and $y=v_{p^{-1}(0)}=v_0$, so $u_0y, u_1x \in E(G)$. Now we describe a $2$--bisection of $G$. Let $c: V(G) \longrightarrow \{B,W\}$ be a $2$--colouring such that:

\begin{itemize}
\item $c(u_0)=c(u_1):=W$;
\item $c(u_i):=W$ if $i$ is odd;
\item $c(u_i):=B$ if $i$ is even ($\ne 0$);
\item $c(x)=c(y):=B$.
\end{itemize}

We need to describe the $2$--colouring $c$ on the rest of the inner cycle in such a way that $c$ is a $2$--bisection.

Since the inner cycle $C_2$ is odd, the vertices $x$ and $y$ divide $C_2$ into two paths, one of even length, say $P_1$, and one of odd length, say $P_2$.
Let $x_0$ and $y_0$ be the neighbours of $x$ and $y$ in $P_1$ and let $x_1$ and $y_1$ be the neighbours of $x$ and $y$ in $P_2$, respectively.

Since $||p_{i+1} - p_{i}||$ is even and $||p^{-1}_{i+1} - p^{-1}_{i}||$ is odd, for all $i$, we have that the vertex $u_h \in C_1$, with $u_hy_1 \in E(G)$, is such that $h$ is even. Then $y_1$ is adjacent to a black vertex of $C_1$. Analogously, the vertex $u_k \in C_1$, with $u_kx_1 \in E(G)$, is such that $k$ is even, hence also $x_1$ is adjacent to a black vertex in $C_1$. Since $x_1$ has two black neighbours ($x$ and $u_k$), we must set $c(x_1):=W$. Similarly, $c(y_1):=W$.

Now consider $x_0u_l, y_0u_m \in E(G)$, with $u_l, u_m \in C_1$. Then, again because $||p_{i+1} - p_{i}||$ even for all $i$, $l$ and $m$ are both odd, which implies that $c(u_l)=c(u_m)=W$ have been set.

Since the number of white vertices in $C_1$ is one more that the number of black vertices, we need the opposite situation in $C_2$. We will obtain this condition in two different ways according to necessity: Type $I$, where $c(x_0)=c(y_0)=W$, there are two consecutive black vertices somewhere in $P_2$ and the rest is alternating; and Type $II$, where $c(x_0)=c(y_0)=B$, there are two consecutive white vertices somewhere in $P_2$ and the rest is alternating, except for the edges $xx_0$ and $yy_0$ which have both ends black (see Figure~\ref{configI&II}).

\begin{figure}[h]
	\centering
		\includegraphics[height=5cm]{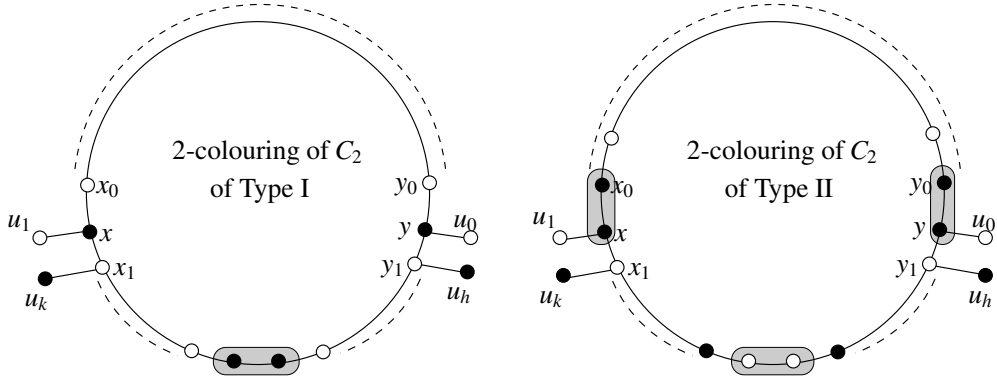}
		\caption{Type I and Type II colouring of $C_2$}\label{configI&II}		
\end{figure}

Since the $C_2(x_1,\widehat{x},y_1)$--path is odd, then there is an even number of vertices between $x_1$ and $y_1$ in this path, say $x_2, \ldots, x_t,y_t, \ldots y_2$.

Now we can complete the $2$--colouring $c$ in order to obtain a $2$--bisection.

{\sc Trivial case:} there is no vertex between $x_1$ and $y_1$. Then $x_1,y_1$ are two consecutive white vertices in $C_2$, so we set $c(x_0)=c(y_0)=B$ and complete $c$ as a Type $II$ colouring.

{\sc General case:} If the two neighbours of $x_2$ and $y_2$ in $C_1$ have both the same colour of $x_1$ and $y_1$, i.e.\ white, then we must set $c(x_2)=c(y_2)=B$ to avoid monochromatic components of order larger than two. Otherwise, we can colour either $c(x_2)=W$ or $c(y_2)=W$, obtaining two consecutive white vertices (either $x_1x_2$ or $y_1y_2$) and we can complete $c$ as a Type $II$ colouring.

We repeat the process by looking at the neighbours of $x_3, y_3$ in $C_1$, and so on. At each step, either the colouring of $x_i$ and $y_i$ is uniquely determined, or we can obtain two consecutive vertices of the same colour. If these two vertices are white, we complete $c$ as a Type $II$ colouring; if these two consecutive vertices are black, we complete $c$ as a Type $I$ colouring.

Finally, if the colouring of all vertices $x_2, \ldots, x_t,y_t, \ldots y_2$ is forced, then $x_t$ and $y_t$ share the same colour; we complete the colouring $c$ as a Type I or a Type II colouring according to the colour of $x_t$ and $y_t$. This completes the proof.
\end{proof}

To conclude this section, we would like to remark that a direct proof of Ando's Conjecture for cycle permutation graphs is much easier. If we do not care about the order of connected components, then we can consider the bisection $B=V(C_1)$ and $W=V(C_2)$ and it is clear that such a bisection has isomorphic parts.

\subsection{1--connected cubic graphs}\label{1connectedBanLinial}

In the previous sections, we have discussed the problem of constructing bisections of cubic graphs with isomorphic parts. In particular, we have observed that a $2$--bisection is always a bisection with isomorphic parts (see Proposition~\ref{BanLinialimpliesAndo}).
As mentioned before, only one example of bridgeless cubic graph without a $2$--bisection is known (i.e.\ the Petersen graph). On the other hand, an infinite family of $1$--connected cubic graphs with no $2$--bisection is constructed in~\cite{EMT15}. In this section, we prove that all members of this family (and also of a slightly more general family) admit a $3$--bisection with isomorphic parts.

\begin{figure}[htb]
\centering \includegraphics[scale=0.8]{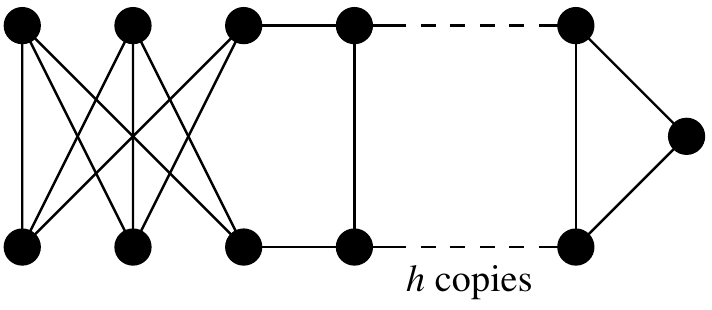} \caption{The module $L_h$}
\label{fig:Lh}
\end{figure}

First of all, we present the family of examples constructed in~\cite{EMT15}:
take the module $L_h$ ($h\ge 0$) depicted in Figure~\ref{fig:Lh}.
Note that $L_0$ is $K_{3,3}$ with on edge subdivided with a new vertex.

Let $T_{ijk}$ be the graph obtained by taking $L_i$, $L_j$ and $L_k$
and adding a new vertex adjacent to the three vertices of degree 2.

It is proved in~\cite{EMT15} that $T_{ijk}$ does not admit a $2$--bisection for all possible non--negative values of $i,j$ and $k$.

\begin{figure}[htb]
\centering
\includegraphics[scale=0.6]{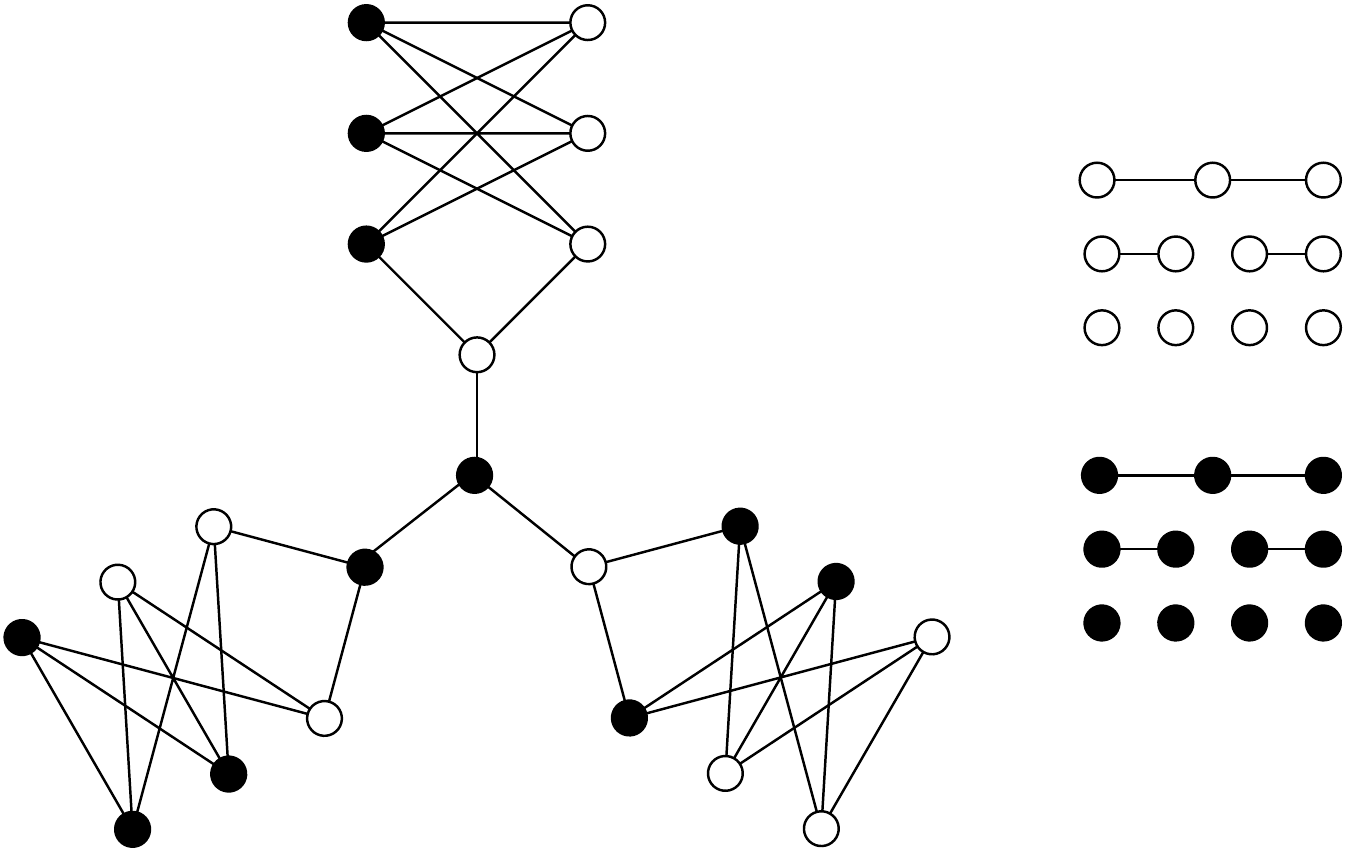} \caption{A $3$--bisection of $T_{000}$ with isomorphic parts}
\label{fig:T000}
\end{figure}

\begin{proposition}\label{notcounterando}
 For any $i,j,k\ge 0$, the graph $T_{ijk}$ admits a $3$--bisection with isomorphic parts.
\end{proposition}

\begin{proof}
The basic step is the $3$--bisection of the graph $T_{000}$ depicted in Figure~\ref{fig:T000}.
The proof proceeds by showing that if we replace the three copies of the module $L_0$ with copies of $L_i$, $L_j$ and $L_k$, for arbitrary non--negative values of $i$,$j$ and $k$, then we can still obtain a $3$--bisection with isomorphic parts.

1) Replace the central copy of $L_0$ in $T_{000}$ with a copy of $L_i$ coloured as in Figure~\ref{fig:Li}.
\begin{figure}[htb]
\centering
\includegraphics[scale=0.6]{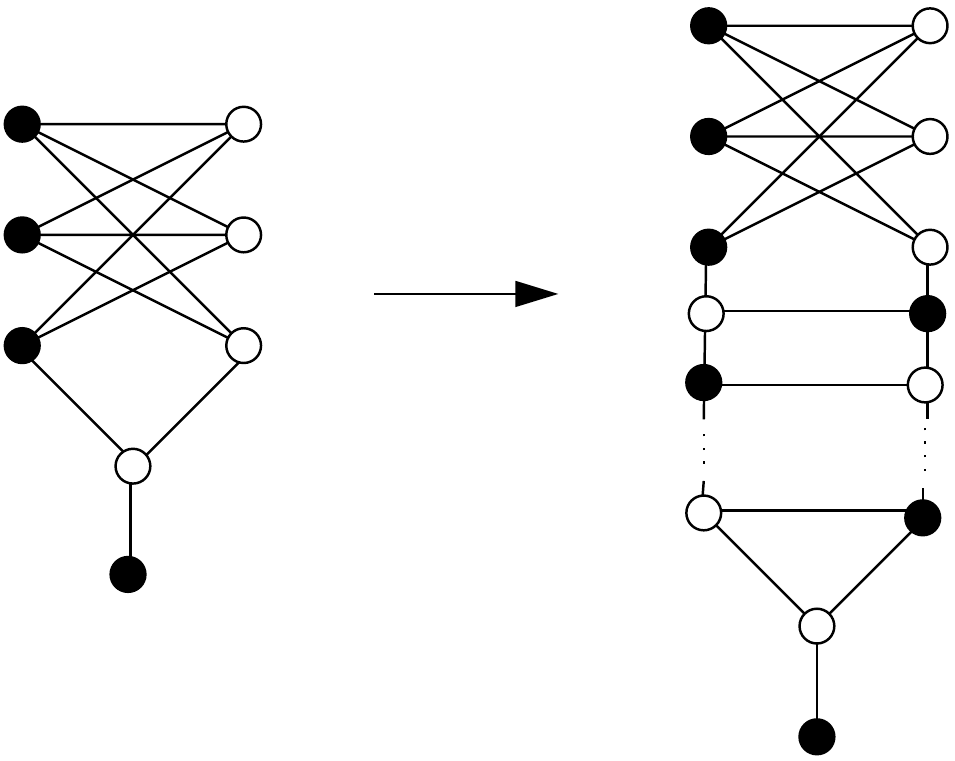} \caption{}
\label{fig:Li}
\end{figure}
The bisection of the obtained graph still has isomorphic parts since we only add isolated vertices in each monochromatic subgraph.

2) Replace the copy of $L_0$ on the right with a copy of $L_j$ coloured as in Figure~\ref{fig:Lj}.
\begin{figure}[htb]
\centering
\includegraphics[scale=0.6]{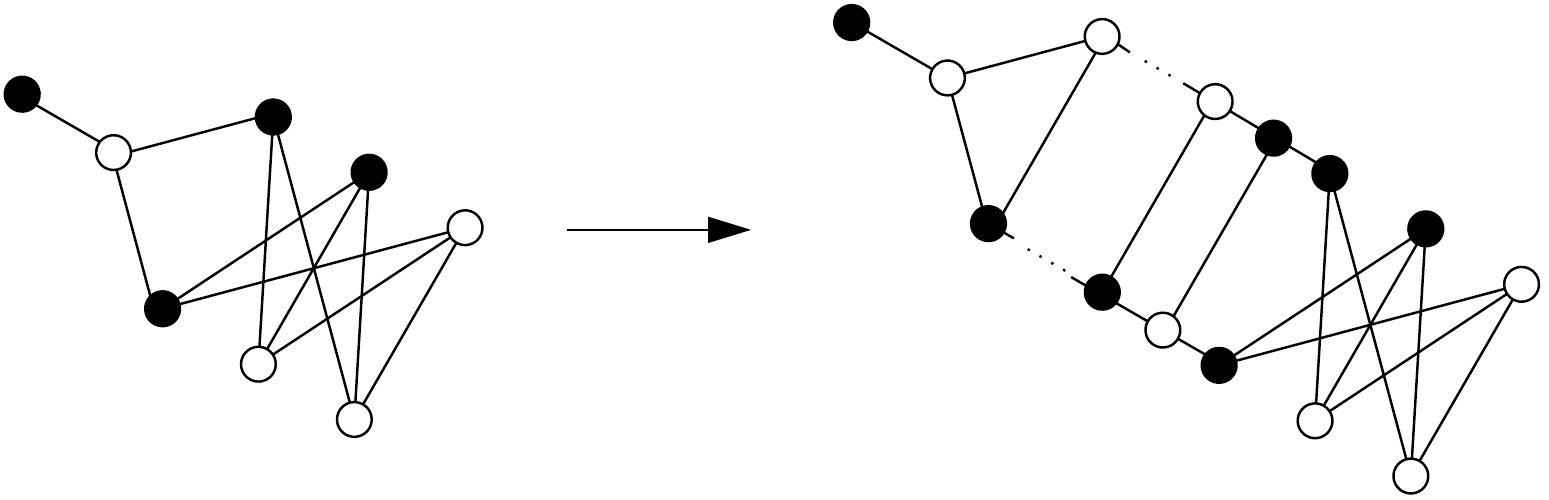} \caption{}
\label{fig:Lj}
\end{figure}
The bisection of the graph obtained still has isomorphic parts since we replace an isolated vertex with a path of length one in each part and we add the same number of isolated vertices in each of them.

3) Replace the copy of $L_0$ on the left with a copy of $L_k$ coloured as in Figure~\ref{fig:Lk} according to the parity of $k$.
\begin{figure}[htb]
\centering
\includegraphics[scale=0.6]{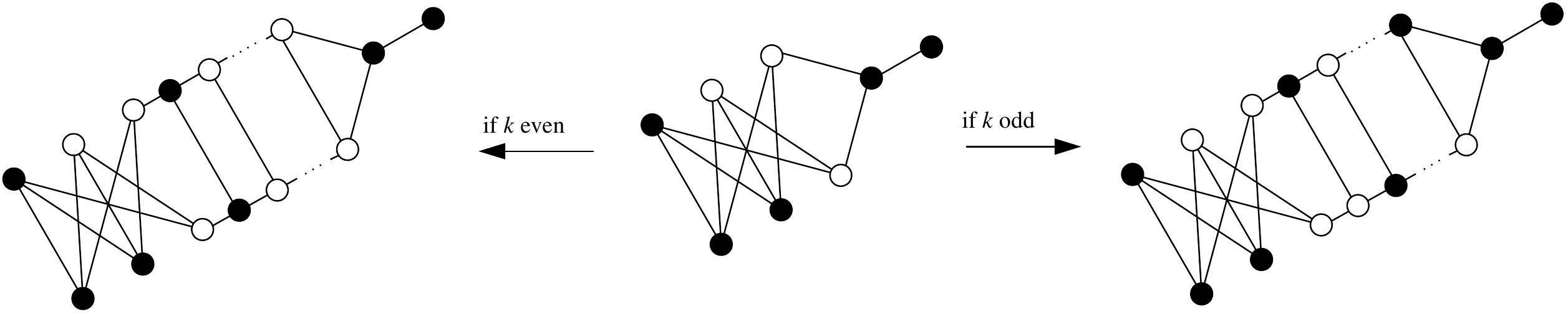} \caption{}
\label{fig:Lk}
\end{figure}
In both cases the bisection of the graph obtained still has isomorphic parts. If $k$ is even, we add the same number of paths of length one in each part, if $k$ is odd we replace a path of length one with a path of length two in each part and we add the same number of isolated vertices in each of them.
\end{proof}

\begin{figure}[htb]
\centering
\includegraphics[scale=0.4]{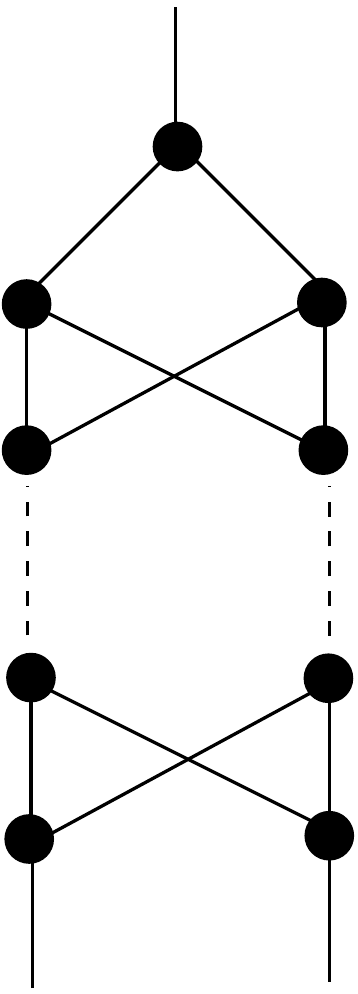} \caption{A gadget to construct new families of cubic graphs without a $2$--bisection}
\label{fig:ker2}
\end{figure}

It is easy to prove that if we replace the central vertex of $T_{ijk}$ with a copy of the graph in Figure~\ref{fig:ker2} having $4h+1$ vertices, we still obtain a cubic graph without a $2$--bisection. So this leads to additional infinite families of cubic graphs without a $2$--bisection.

\begin{figure}[htb]
\centering
\includegraphics[scale=0.4]{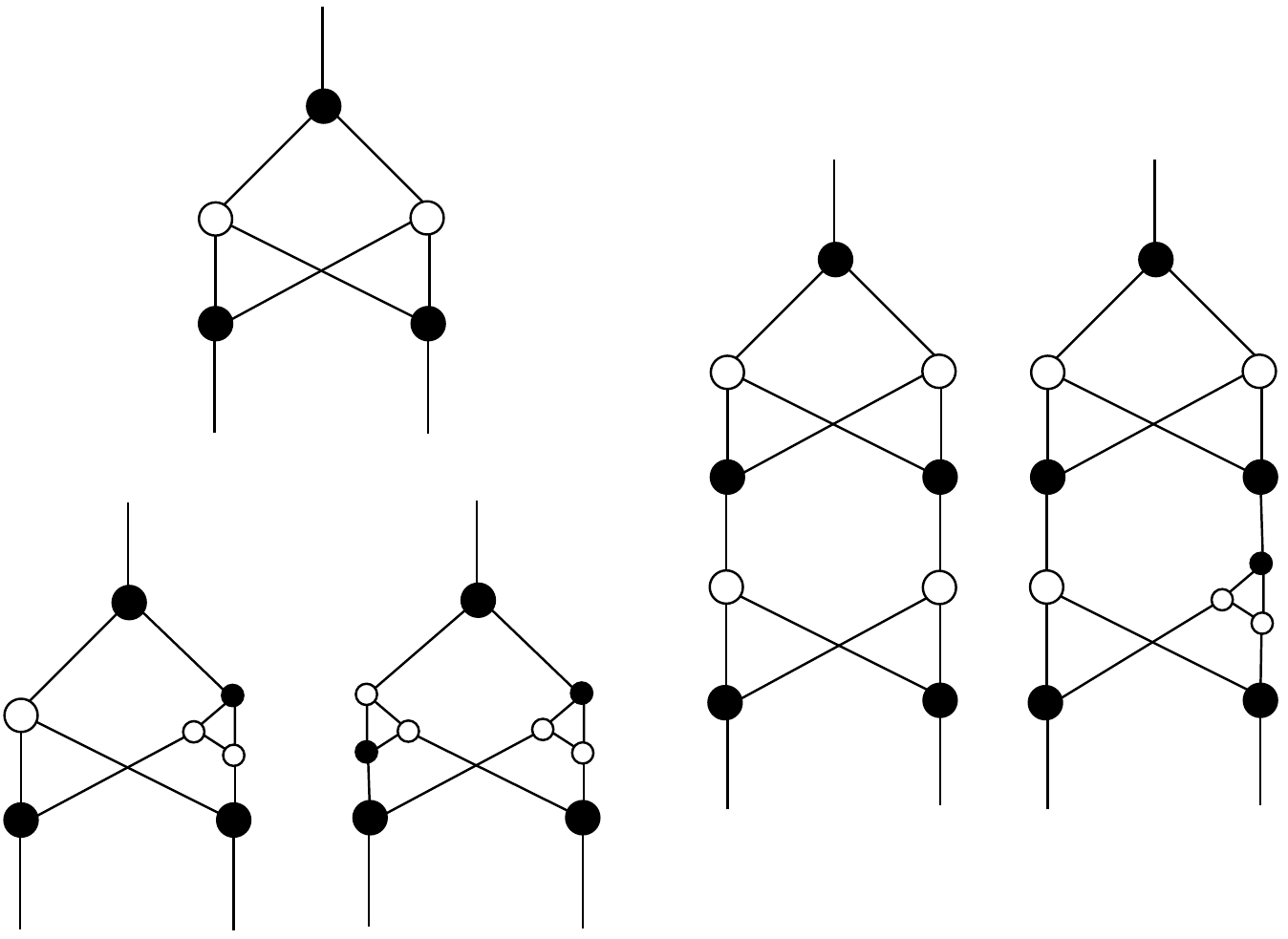} \caption{All cubic graphs of order at most 32 (except for the Petersen graph) with no $2$--bisection can be obtained by replacing the central vertex in $T_{ijk}$ with one of these subgraphs }
\label{fig:ker}
\end{figure}

Moreover, all $1$--connected graphs of order at most 32 without a $2$--bisection are computationally determined in Section~\ref{CompRes}. We verified that all of them are obtained by replacing the central vertex of $T_{ijk}$ with a copy of one of the graphs from Figure~\ref{fig:ker}.
Note that all of the central parts from Figure~\ref{fig:ker} are either a member of the family of subgraphs shown in Figure~\ref{fig:ker2} or are obtained by a member of such a family by replacing some of the vertices with triangles. We cannot say in general if there exist any larger examples which do not fit into such a description.

By colouring the vertices of these subgraphs as shown in Figure~\ref{fig:ker}, one can see that Proposition~\ref{notcounterando} can be extended to prove that each member of these additional families admits a $3$--bisection with isomorphic parts. Hence, none of them is a counterexample for Ando's Conjecture.

Finally, let us remark that none of the graphs presented in this section without a $2$-bisection admits a perfect matching. Thus, none of them is a counterexample to the following stronger version of Ban-Linial Conjecture proposed in~\cite{EMT15}.

\begin{conjecture}[\cite{EMT15}]
The Petersen graph is the only cubic graph that admits a perfect matching, but no $2$-bisection.
\end{conjecture}

\subsection{Strong Ando Conjecture}\label{subsect:strongando}

We conclude this section with some properties of bisections with isomorphic parts that will allow us to state a strengthening of Ando's conjecture that will be discussed in Section~\ref{strongwormaldconj}.

Recall that, as was remarked at the beginning of this section, in all $2$--bisections the two isomorphic parts are obviously linear forests, since the connected components are either isolated vertices or isolated edges. Moreover, we have verified that all known graphs without a $2$--bisection admit a $3$--bisection with isomorphic parts which are linear forests (this is not trivial since a connected component of order $3$ could be a $3$--cycle) except the Petersen graph (see Proposition~\ref{notcounterando} and Observation~\ref{obs:bisectionwithshortlinearforests}). Finally, note that the Petersen graph admits a $4$--bisection with isomorphic parts which are linear forests (see Figure~\ref{figure1}(c)).

In the following remark, we summarise all of these facts:

\begin{remark}\label{bisectionandisoparts}

\begin{enumerate}
\item[(i)] In all $2$--bisections the two monochromatic induced subgraphs are isomorphic and they are obviously linear forests;

\item[(ii)] All known graphs without a $2$--bisection admit a $3$--bisection with isomorphic parts which are linear forests, except the Petersen graph;

\item[(iii)] The Petersen graph admits a $4$--bisection with isomorphic parts which are linear forests.
\end{enumerate}

\end{remark}

This gives rise to the following natural question: {\em does every cubic graph admit a bisection with isomorphic parts and with the additional property that these parts are linear forests?}

To approach this problem, first of all, we observe that we can assume that both parts have maximum degree two.

\begin{proposition}\label{bisectionwithisopartsdegree2}
Let $G$ be a cubic graph admitting a bisection with isomorphic parts. Then $G$ admits a bisection with isomorphic parts having maximum degree $2$.
\end{proposition}

\begin{proof}
Take a bisection with isomorphic parts and with the minimum number of vertices of degree $3$ in each induced monochromatic subgraph.
Let $v_1$ be a vertex of degree $3$ in the subgraph induced by the black vertices, and let $v_2$ be the corresponding vertex in the isomorphic monochromatic white subgraph. Switching colours between $v_1$ and $v_2$ would decrease the number of such vertice, contradicting the minimality assumption.
\end{proof}

Proposition~\ref{bisectionwithisopartsdegree2} and Remark~\ref{bisectionandisoparts}  lead us to state the following strong version of Ando's Conjecture:

\begin{conjecture}[Strong Ando Conjecture]\label{strongandoconj}
Every cubic graph admits a bisection such that the two induced subgraphs are isomorphic linear forests.
\end{conjecture}

We are curious about the optimal value of $k$ in the previous conjecture and we leave it as an open problem. At this stage, we only remark that $k=3$ is always sufficient in all known examples, except for the Petersen graph where we need $k=4$.

\begin{problem}
Does every cubic graph admit a $4$--bisection such that the two induced subgraphs are isomorphic linear forests?
\end{problem}

\section{Strong Wormald (edge) colourings vs Strong Ando (vertex) colourings}\label{strongwormaldconj}

In this section, as mentioned in the Introduction,  we develop a second and different approach introducing the concept of {\em (Strong) Wormald Colouring},  to present further evidence in support of Ando's Conjecture (i.e.\ Conjecture~\ref{ando}), and its stronger version (Conjecture~\ref{strongandoconj_colouringversion}).

In the Introduction, we have stated a Conjecture by Wormald (cf.\ Conjecture~\ref{Wconj}) in terms of isomorphic linear forests (i.e.\ {\em let $G$ be a cubic graph with $|V(G)|$ $\equiv 0 \pmod 4$. Then, there exists a linear partition of $G$ into two isomorphic linear forests}).
Note that, in other words, Conjecture~\ref{Wconj} states that {\em there exists a $2$--edge colouring of $G$, $c_E : E(G) \longrightarrow \{B,W\}$ such that both monochromatic subgraphs are linear forests and they are isomorphic}. We define such a colouring $c_E$  to be a \emph{Wormald colouring}.

Analogously, we can restate Ando's Conjecture in terms of colourings as follows:

\begin{conjecture}[Ando's conjecture~\cite{A17} (Colouring Version)]\label{andoreformulated}
A cubic graph $G$ admits a $2$--vertex colouring of $G$,
$c_V : V(G) \longrightarrow \{B,W\}$, such that the monochromatic induced subgraphs are isomorphic.
\end{conjecture}

From now on, a $2$--colouring of $G$, $c_V : V(G) \longrightarrow \{B,W\}$, such that the monochromatic induced subgraphs are isomorphic will be called an \emph{Ando colouring}. Trivially, in order to have isomorphic monochromatic induced subgraphs such a $2$--colouring is a bisection.

It seems immediately clear that there should be a link between a Wormald Colouring of the edges of a cubic graph $G$ and an Ando Colouring of the vertices of $G$. In order to establish a strong connection between the two problems, we need to strengthen some of the previous definitions. In particular, we need to exclude induced paths of length one (i.e. a path with one edge and two vertices) in a Wormald Colouring.

\begin{definition}
Let $G$ be a cubic graph with $|V(G)|$ $\equiv 0 \pmod 4$. A Strong Wormald Colouring is an edge colouring $c_E : E(G) \longrightarrow \{B,W\}$ of $G$ such that the monochromatic induced subgraphs are isomorphic linear forests with paths of length at least $2$.
\end{definition}

Finally, we define a Strong Ando Colouring as an Ando Colouring such that the monochromatic induced subgraphs are not only isomorphic but also linear forests. More precisely:

\begin{definition}
Let $G$ be a cubic graph. A Strong Ando Colouring is a vertex colouring $c_V : V(G) \longrightarrow \{B,W\}$ of $G$ such that the monochromatic induced subgraphs are isomorphic linear forests.
\end{definition}

Obviously, a Strong Ando Colouring is an Ando Colouring. Moreover, we can restate Conjecture~\ref{strongandoconj} in terms of colouring as follows:

\begin{conjecture}[Strong Ando Conjecture (Colouring Version)]\label{strongandoconj_colouringversion}
Every cubic graph admits a Strong Ando Colouring.
\end{conjecture}

The idea of introducing the concept of a Strong Wormald Colouring is mainly due to the following proposition.

\begin{proposition}\label{swortosandocol}
Let $G$ be a cubic graph with $|V(G)|$ $\equiv 0 \pmod 4$ admitting a Strong Wormald Colouring. Then, $G$ admits a Strong Ando Colouring (and thus an Ando Colouring).
\end{proposition}

\begin{proof}
Let $c_E : E(G) \longrightarrow \{B,W\}$ be a Strong Wormald Colouring. We construct a Strong Ando Colouring as follows:

$$
c_V(v) := \left\{ \begin{array}{ccc}
                    B & if & v \mbox{ is incident with two black edges of } c_E \\
                    W & if & v \mbox{ is incident with two white edges of } c_E
                  \end{array}
\right.
$$
The vertex colouring $c_V(v)$ is indeed a Strong Ando Colouring: every connected component of a monochromatic induced subgraph is a path, since it is obtained starting from a path of the induced linear forests of the Strong Wormald Colouring by removing the two final edges. Hence, the two monochromatic induced subgraphs are linear forests, and the previous argument implicitly proves that they are also isomorphic.
\end{proof}

As far as we know, Wormald's Conjecture (i.e.\ every cubic graph of order  $0 \pmod 4$ has a Wormald Colouring) is open  and we have no counterexamples for the Strong Ando Conjecture (i.e.\ every cubic graph has a Strong Ando Colouring, and so also an Ando Colouring). What happens if we consider Strong Wormald Colourings instead?
Is it true that every cubic graphs of order congruent to $0 \pmod 4$ has a Strong Wormald Colouring? The answer is negative in general (see Observation~\ref{obs:strong_wormald} for more details). In the Appendix we give a description of all $2$-connected examples of order at most $28$ and we construct an infinite family of $2$--connected graphs without a Strong Wormald Colouring. The family of these examples seems to represent a very thin subclass of the class of all cubic graphs and, even if they do not admit a Strong Wormald Colouring, all of them admit a $2$--bisection and hence a Strong Ando Colouring.

In order to have an efficient procedure to obtain graphs without a Strong Wormald Colouring, we first characterise bisections which arise from a Strong Wormald Colouring.

Let $G$ be a cubic graph such that $|V(G)|$ $\equiv 0 \pmod 4$. Let $c_E: E(G) \longrightarrow \{B,W\}$ be a Strong Wormald Colouring. Let $c_V : V(G) \longrightarrow \{B,W\}$ be the vertex colouring induced by $c_E$ as described in the proof of Proposition~\ref{swortosandocol}.
From now on we denote by $G'$ the spanning subgraph of $G$ obtained by removing all edges with two ends of the same colour in $c_V$. In other words, $G' = (V(G'), E(G'))$ where $V(G')=V(G)$ and $E(G') = \{ xy \in E(G) : c_V(x) \neq c_V(y) \} = \{ xy \in E(G): \mbox{ either $x$ or $y$ is a leaf in a monochromatic path of } \, c_E\}$.

In the following lemma we describe the connected components of the subgraph $G'$.

\begin{lemma}\label{treeplusedge}
Let $G$ be a cubic graph with $|V(G)|$ $\equiv 0 \pmod 4$ admitting a Strong Wormald Colouring. Each connected component of the subgraph $G'$ of $G$ is bipartite and isomorphic to a tree plus an edge.
\end{lemma}

\begin{proof} Let $c_V : V(G) \longrightarrow \{B,W\}$ be the vertex colouring induced by the Strong Wormald Colouring of $G$ (see the proof of Proposition~\ref{swortosandocol}).
Since all edges with ends of the same colour in $c_V$ were removed from $G$, the graph $G'$ is obviously bipartite with bipartition $c_V$.
Consider the following orientation of $G'$: $e \in E(G')$ is oriented away from the vertex having the same colour in $c_V$ as $e$ in $c_E$.
This orientation cannot produce a vertex of in--degree larger than $1$ in $G'$: take a vertex $v$, the colour in $c_E$ of an inner edge in $v$ is opposite to the colour of $v$ in $c_V$, and we can have only one edge of such a colour incident to $v$ by definition of $c_V$.
Since we do not have in--degree $2$ at any vertex,  every path of $G'$ is a directed path and every cycle $C$ in $G'$ is a directed cycle (and even, since $G'$ is bipartite).
Furthermore, every edge incident to a vertex of $C$, and not in $C$, is directed away from $C$.
Then, every path leaving $C$ cannot come back to $C$.
For the same reason, we cannot have two cycles in the same connected component of $G'$.
Thus, we have proved that if a connected component of $G'$ has a cycle, it is the unique cycle in such a component, in other words it is a tree plus an edge.

Finally, we observe that a connected component cannot be acyclic (i.e. a tree). Each edge incident to a leaf is oriented towards its leaf by definition of $c_V$ and $c_E$, so a path which contains two leaves cannot be a directed path, a contradiction.
Hence, a connected component of $G'$ is bipartite and it contains exactly one cycle.
\end{proof}

Now we prove that every bisection having the properties described in the previous propositions arises from a Strong Wormald Colouring.

\begin{theorem}\label{equivalencebisectionstrongwormald}
Let $G$ be a cubic graph with $|V(G)|$ $\equiv 0 \pmod 4$. The graph $G$ admits a Strong Wormald Colouring if and only if $G$ admits a Strong Ando Colouring such that each connected component of the bipartite subgraph $G'$ of $G$, obtained by removing all edges of $G$ with ends of the same colour, is isomorphic to a tree plus an edge.
\end{theorem}

\begin{proof} One direction is exactly the statement of Lemma~\ref{treeplusedge} and now we prove the other direction.
Let $c_V$ be a Strong Ando Colouring of $G$ with the additional property that each connected component in the bipartite subgraph $G'$ of $G$, obtained by removing all edges of $G$ with ends of the same colour, is isomorphic to a tree plus an edge.

We construct an edge colouring $c_E$ of $G$ starting from $c_V$ and then we prove that $c_E$ is a Strong Wormald Colouring.\\
For every edge $xy$ having both ends of the same colour in $c_V$,  we set
$$ c_E(xy) := \left\{
\begin{array}{ccc}
  B & if & c_V(x)=c_V(y)=B \\
  W & if & c_V(x)=c_V(y)=W
\end{array}
\right.$$

Now we must define $c_E$ on the remaining edges of $G$, which are exactly the edges of $G'$.
Consider a connected component of $G'$: it has a unique even cycle $C$, so colour the edges of $C$ alternately black and white. Now colour the remaining edges following the rule that a black vertex must have two black incident edges (and the same for white vertices). The colouring $c_E$ is uniquely determined up to the colour of each even cycle in $G'$. (Then we have exponentially many different edge colourings in terms of the number of even cycles in $G'$).

Finally, we have to prove that $c_E$ is indeed a Strong Wormald Colouring.

$$
E(G) = E(G') \cup E(G'')
$$

where $E(G')$ are the edges of $G$ with ends of opposite colours in $c_V$, and $E(G'')$ are the edges with ends of the same colour in $c_V$.
The graph $G''$ is the union of two isomorphic linear forests $F_B$ and $F_W$, since $c_V$ is a Strong Ando Colouring.

Now we show how a path $(x_0,x_1,\ldots,x_l)$ of length $\ell \geq 0$ in $F_B$ (and an analogous result holds for $F_W$) corresponds to a path of length $\ell +2$ in $G$.

All vertices $x_i : i=0,\ldots, l $ are black in $c_V$, since they belong to $F_B$. Moreover, their neighbours in $G$ are all white vertices in $c_V$ by definition of $F_B$. Possible inner vertices of the path are all vertices of degree one in $G'$ and so, by definition of $c_E$ all edges incident with the path in these vertices are white. Moreover, the two ends of the path, $x_0$ and $x_l$, are vertices of degree two in $G'$. Exactly one of the two edges incident to $x_0$ (resp. $x_l$) in $G'$ is black in $c_E$, say $e_0$ (resp. $e_l$). The edge $e_0$ (resp. $e_l$) is an isolated edge in $G'$ because its ends distinct from $x_0$ (resp. $x_l$) is white.
Hence, $c_E$ induces linear forests which are obtained from the linear forests induced by  $c_V$ adding two edges at the ends of each path. This implies that the two linear forests are isomorphic and have length at least $2$.

Therefore, $c_E$ is a Strong Wormald Colouring and this proves the theorem.
\end{proof}

Note that Theorem~\ref{equivalencebisectionstrongwormald} not only has relevance from a theoretical point of view, but as will be explained in Section~\ref{comput:wormald}, it also permits to design a more efficient algorithm to search for Strong Wormald Colourings.

\begin{figure}[h]
\centering
\includegraphics[width=9cm]{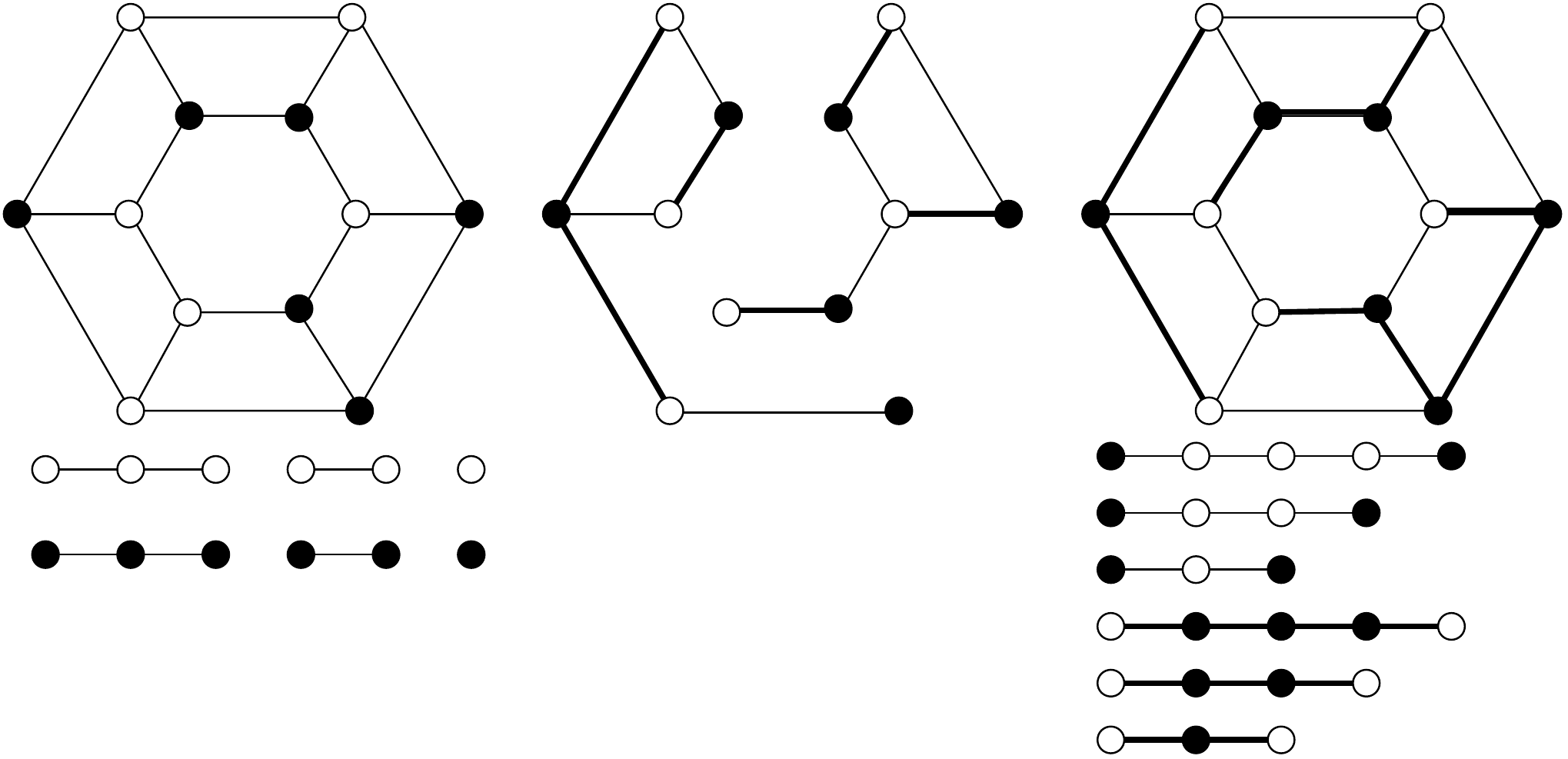}
\caption{The reconstruction of a Strong Wormald Colouring obtained following the proof of Theorem~\ref{equivalencebisectionstrongwormald}}\label{fig:prism12vertices}
\end{figure}

We conclude this section with an example of a reconstruction of a Strong Wormald Colouring of the prism with 12 vertices starting from a suitable Strong Ando Colouring. Figure~\ref{fig:prism12vertices} shows a Strong Ando Colouring on the left, an edge colouring of the bipartite subgraph $G'$ with all connected components isomorphic to a tree plus an edge in the center, and the Strong Wormald Colouring on the right.

\subsection{Strong Wormald Colouring of cubic graphs: $|V| \equiv 2 \pmod 4$}

The discussion of previous section creates an evident connection between Wormald's Conjecture and (Strong) Ando's Conjecture, but it also leaves a large gap between them. Indeed, Wormald's Conjecture is only for cubic graphs of order a multiple of four while (Strong) Ando's Conjecture is stated for all cubic graphs. The purpose of this section is to fill this gap by defining a Strong Wormald Colouring for cubic graphs of order  $2 \pmod 4$.

The original Wormald Conjecture is not stated  for cubic graphs of order not a multiple of four since they have an odd number of edges, thus the two monochromatic parts cannot have the same number of edges and hence they cannot be isomorphic.

We would like to give a definition of Strong Wormald Colouring for cubic graphs with an odd number of edges such that the existence of such a colouring  implies the existence of a Strong Ando Colouring as we did previously. Moreover, we would like to have a definition that is satisfied by a large fraction of cubic graphs as it happens in the case of an even number of edges.

Our very simple idea is to consider an edge colouring of all edges of the graph except one, with an additional property on the uncoloured edge. The presence of this additional property will be exactly what we need to construct a Strong Ando Colouring starting from this edge colouring.

\begin{definition}\label{strongandocolouring}
Let $G$ be a cubic graph of order $2 \pmod 4$ and let $xy \in E(G)$ be an edge of $G$. A Strong Wormald Colouring of $G$ is an edge colouring $c_E : E(G) - \{xy\} \longrightarrow \{B,W\}$ of $G$ such that the monochromatic induced subgraphs are isomorphic linear forests with paths of length at least $2$, and, moreover, the vertex $x$ is incident to two white edges and the vertex $y$ is incident to two black edges.
\end{definition}

This definition of Strong Wormald Colouring allows to state the analogous of Proposition~\ref{swortosandocol}.
We omit the proof since it proceeds exactly like the proof of Proposition~\ref{swortosandocol} for order $0 \pmod 4$.

\begin{proposition}\label{swortosandocol2}
Let $G$ be a cubic graph with $|V(G)|$ $\equiv 2 \pmod 4$ admitting a Strong Wormald Colouring. Then, $G$ admits a Strong Ando Colouring (and so an Ando Colouring).
\end{proposition}

Moreover, we have also obtained the following characterization which is similar to the case $ 0 \pmod 4$ of bisections arising from a Strong Wormald Colouring (cf.\ Theorem~\ref{equivalencebisectionstrongwormald}).

\begin{theorem}\label{equivalencebisectionstrongwormald2}
Let $G$ be a cubic graph  with $|V(G)|$ $\equiv 2 \pmod 4$. The graph $G$ admits a Strong Wormald Colouring if and only if $G$ admits a Strong Ando Colouring such that each connected component in the bipartite subgraph $G'$ of $G$, obtained by removing all edges of $G$ with ends of the same colour, is isomorphic to a tree plus an edge except for one that is a tree.
\end{theorem}

The proof is analogous to the proof of Theorem~\ref{equivalencebisectionstrongwormald}. The unique difference is the connected component which is isomorphic to a tree, but this fact easily follows due to the presence of the uncoloured edge.

So far, we have followed exactly the same idea of the case of order $0 \pmod 4$: we have defined a Strong Wormald Colouring, we have proved that the existence of such an edge colouring implies the existence of a (Strong) Ando Colouring and we have shown a characterization in terms of bisections of the graphs admitting a Strong Wormald Colouring.

Note that a large amount of cubic graphs seem to have a Strong Wormald Colouring. Indeed,
the possibility of choosing the edge $xy$ produces a better result than in the case $0 \pmod 4$. In particular, using a computer program, we have checked all cubic graphs having at most $26$ vertices of order $2 \pmod 4$ and {\em all} of them admit a Strong Wormald Colouring (see Section~\ref{comput:wormald} for details).

\begin{remark}
As suggested by one of the anonymous referees, we may apply a similar approach to the case $0 \pmod 4$ by removing a pair of edges.
In other words, we may look for two edges $e$ and $f$ of $G$ such that there exist a Strong Wormald Colouring of $G-\{e,f\}$ and the ends of $e$ and $f$ satisfy a condition analogue to the one in Definition~\ref{strongandocolouring}. The existence of such a pair in $G$ permit to prove that $G$ admits a Strong Ando Colouring, exactly in the same way we did it in Proposition~\ref{swortosandocol2}.

A computational check of all cubic graphs of order up to 24 shows that all of them admit such a pair except for two instances. The first one is a graph of order $12$ obtained by replacing the three vertices of $K_{3,3}$ in the same bipartition class with triangles. The second one has order $16$ and is the smallest simple cubic graph with no perfect matching. These computations seem to suggest a new possible approach to prove the Strong Ando Conjecture.
Unfortunately, at the moment, we could not prove in general the existence of such a pair of edges in every cubic graph of order $0 \pmod 4$.
\end{remark}


We leave the existence of a cubic graph of order $2 \pmod 4$ without a Strong Wormald Colouring as an open problem.

\begin{problem}\label{prob:strenghtstrongwormald}
Does every cubic graph of order $2 \pmod 4$ have a Strong Wormald Colouring?
\end{problem}

Note that, by Proposition~\ref{swortosandocol2}, a positive answer to Problem~\ref{prob:strenghtstrongwormald} would imply (Strong) Ando's conjecture for all cubic graphs of order order $2 \pmod 4$.

\section{Computational Results}\label{CompRes}

In this section we present our computational results. These results also give evidence to support the correctness of the conjectures which were introduced in the previous sections. More specifically, in Section~\ref{comput:bl} we present our computational results on the conjecuture of Ban and Linial,  in Section~\ref{comput:ando} on Ando's conjecture and finally in Section~\ref{comput:wormald} on Wormald's conjecture.

We wrote programs for testing each of these conjectures or open problems. The source code of these programs can be downloaded from~\cite{bisection-site}.

In each case, we used the generator for cubic graphs \textit{snarkhunter}~\cite{snarkhunter-site,brinkmann_11} to generate all cubic graphs up to a given number of vertices and then applied our programs for testing the conjectures on the generated graphs. In each case we went as far as computationally feasible.

\subsection{Computational results on Ban--Linial's Conjecture}
\label{comput:bl}

We wrote a program for testing if a given graph has a 2--bisection. We generated all cubic graphs up to 32 vertices and then tested which of these graphs do not admit a $2$--bisection. This yielded the following observation.

\begin{observation} \label{obs:2bisection_cubic}
There are exactly 34 graphs among the cubic graphs with at most 32 vertices which do not admit a $2$--bisection. All of these graphs, except the Petersen graph, have connectivity $1$.
\end{observation}

So this implies that the smallest counterexample to the Ban--Linial Conjecture must have at least 34 vertices.
The counts of the graphs from Observation~\ref{obs:2bisection_cubic} can be found in Table~\ref{table:num_no_2bisection}. We verified that all of them are obtained by replacing the central vertex of $T_{ijk}$ with a copy of one of the graphs from Figure~\ref{fig:ker} (see Section~\ref{1connectedBanLinial} for more details).

The graphs from Observation~\ref{obs:2bisection_cubic} can also be downloaded and inspected in the database of interesting graphs from the \textit{House of Graphs}~\cite{BCGM} by searching for the keywords ``no 2--bisection''.

\begin{table}[ht!]
\centering
		\begin{tabular}{| c | c |}
		\hline
		Order & \# no $2$--bisection\\
		\hline
$0-8$	& 	0\\
10  &   1\\
$12-20$ & 0\\
22  &   1\\
24  &   1\\
26  &   3\\
28  &   5\\
30  &   9\\
32  &   14\\
		\hline
		\end{tabular}
		\caption{Number of cubic graphs which do not admit a $2$--bisection.}
		\label{table:num_no_2bisection}
\end{table}

Recall from Theorem~\ref{BLthm} that Ban and Linial have proved that every $3$--edge colourable cubic graph admits a $2$--bisection. Snarks (as mentioned earlier) are an important subclass of bridgeless not $3$--edge colourable cubic graphs and in~\cite{BGHM} Brinkmann, H{\"a}gglund, Markstr{\"o}m and the second author have determined all snarks up to $36$ vertices.

We now also tested which snarks of order 34 and 36 admit a $2$--bisection, which led to the following observation.

\begin{observation}\label{obs:snark_2bisect}
The Petersen graph is the only snark up to 36 vertices which does not admit a $2$--bisection.
\end{observation}

So this provides further evidence to support the correctness of Conjecture~\ref{banlinial}. Note that it follows from Theorem~\ref{thm:circflownr} that it is in fact sufficient only to test snarks with circular flow number $5$ for Observation~\ref{obs:snark_2bisect}. Together with Mattiolo, the second and fourth author have recently determined all snarks with circular flow number $5$ up to 36 vertices in~\cite{GMM}. (Though note that there is no evidence that a minimal counterexample to the Ban--Linial Conjecture would be cyclically $4$-edge-connected or have girth at least $5$).

\subsection{Computational results on Ando's Conjecture}
\label{comput:ando}

Recall that an easy counting argument shows that every $2$--bisection is a bisection with isomorphic parts. So a possible counterexample to Ando's conjecture must be a graph that does not have a $2$--bisection.

We extended our program to test (the strong version of) Ando's conjecture and applied it to the list of all cubic graphs up to 32 vertices which do not have a $2$--bisection from Observation~\ref{obs:2bisection_cubic}. This led to the following observation.

\begin{observation}\label{obs:bisectionwithshortlinearforests}
All 34 graphs of Observation~\ref{obs:2bisection_cubic} (i.e.\ all cubic graphs which do not admit a $2$--bisection up to 32 vertices), except the Petersen graph, have a $3$--bisection such that the two isomorphic induced monochromatic graphs are linear forests. The Petersen graph has no such bisection, but does admit a $4$--bisection such that the two isomorphic induced monochromatic graphs are linear forests.
\end{observation}

So this implies the following:

\begin{corollary}
The Strong version of Ando's conjecture (i.e.\ Conjecture~\ref{strongandoconj}) (and thus also Ando's original conjecture) does not have any counterexamples with less than 34 vertices.
\end{corollary}

\subsection{Computational results on Wormald's Conjecture}
\label{comput:wormald}

In this section we describe our algorithm to test if a cubic graph admits a \textit{Strong Wormald Colouring} and our results obtained by applying it to the complete lists of cubic graphs.

We implemented an algorithm which tries to find a Strong Wormald Colouring directly by edge colouring the graphs (see later), but it turns out that an algorithm which constructs bisections based on Theorem~\ref{equivalencebisectionstrongwormald} and Theorem~\ref{equivalencebisectionstrongwormald2} is more efficient. The latter algorithm works as follows.

For every cubic graph $G$, we compute all bisections of $G$ until a bisection is found which admits a Strong Wormald Colouring. A bisection leads to a Strong Wormald Colouring if and only if the following conditions (1),(3),(4) are fulfilled.
We add a further condition, i.e.\ Condition (2), which is obviously weaker than Condition (3), but it is added for efficiency reasons.
In particular, the four conditions are listed and executed in this order to obtain a better computational performance.

\begin{enumerate}
\item Both monochromatic induced subgraphs are linear forests;
\item Let denote by $(l_0,\ldots,l_t)$ the vector of lengths of the paths in one of the monochromatic subgraphs. Then for both monochromatic subgraphs, the following must hold: $\sum\limits_{i} (l_i+2) = \lfloor \frac{|E(G)|}{2} \rfloor$.
\item The two induced monochromatic subgraphs are isomorphic;
\item Let $G'$ be obtained from $G$ by the removal of all monochromatic edges. If $|V(G)| \equiv 0 \pmod 4$ then in every connected component of $G'$ the number of edges equals the number of vertices ( and the component is therefore a tree plus an extra edge). If $|V(G)| \equiv 2 \pmod 4$ the above holds for every component of $G'$, except for one component which is a tree.


\end{enumerate}

We tested which of the cubic graphs of order $0 \pmod 4$ do not have a Strong Wormald Colouring. This led to the following observation.

\begin{observation}\label{obs:strong_wormald}
There are exactly 131 graphs without a Strong Wormald Colouring among the cubic graphs of order $0 \pmod 4$ and at most 28 vertices.
\end{observation}

The graph counts and the connectivity statistics of the graphs from Observation~\ref{obs:strong_wormald} can be found in Table~\ref{table:counterex_strong_wormald}. All $2$--connected examples are presented in Figure~\ref{fig:2connectedwithoutstrongwormald} from the Appendix (where we also construct an infinite family of $2$--connected graphs without a Strong Wormald Colouring).
The graphs from Observation~\ref{obs:strong_wormald} can also be downloaded and inspected in the database of interesting graphs from the \textit{House of Graphs}~\cite{BCGM} by searching for the keywords ``no strong Wormald colouring''.

\begin{table}[ht!]
\centering
		\begin{tabular}{| c | c | c | c | c |}
		\hline
		Order & Conn.\ 1 & Conn.\ 2 & Conn.\ 3 & Total\\
		\hline
4  &  0  &  0  &  0  &  0 \\
8  &  0  &  0  &  0  &  0 \\
12  &  0  &  0  &  0  &  0 \\
16  &  1  &  1  &  1  &  3 \\
20  &  20  &  3  &  1  &  24 \\
24  &  18  &  1  &  0  &  19 \\
28  &  72  &  12  &  1  &  85 \\
		\hline
		\end{tabular}
		\caption{Number of cubic graphs of order $0 \pmod 4$ up to 28 vertices which  do not admit a strong Wormald Colouring. Conn.\ $k$ lists the counts of the graphs with connectivity~$k$.}
		\label{table:counterex_strong_wormald}
\end{table}

We also implemented a second independent algorithm to test if a given graph has a Strong Wormald Colouring by edge colouring the graph. We also executed this edge colouring algorithm on all cubic graphs with $0 \pmod 4$ vertices up to 24 vertices and this yielded exactly the same results as in Table~\ref{table:counterex_strong_wormald}.

We also verified that the 131 graphs from Observation~\ref{obs:strong_wormald} without a Strong Wormald Colouring are not counterexamples to the original conjecture of Wormald, that is: they admit a Wormald Colouring. So we can state the following:

\begin{corollary}
Wormald's Conjecture (i.e.\ Conjecture~\ref{Wconj}) does not have any counterexamples with less than 32 vertices.
\end{corollary}

We also investigated Strong Wormald Colourings for graphs of order $2 \pmod 4$. This led to the following observation.

\begin{observation}\label{obs:strong_wormald_2mod4}
Every cubic graph of order $2 \pmod 4$ and at most 26 vertices admits a Strong Wormald Colouring (and hence a (Strong) Ando Colouring).
\end{observation}

Also here the results from Observation~\ref{obs:strong_wormald_2mod4} were independently obtained and confirmed by the bisections algorithm described above and by the algorithm which edge colours the graph.

\section{A related problem on linear arboricity}\label{thomassenprob}

As a by--product of our study of $2$--edge colourings of cubic graphs having linear forests as monochromatic components,  we obtain a negative answer to a problem about linear arboricity posed by Jackson and Wormald in~\cite{JakWor} (cf.\ also~\cite{Tho}). Furthermore, we also pose new problems and conjectures on this topic.

The $k$--linear arboricity of a graph $G$, introduced by Habib and P\'eroche~\cite{HabPer}, is the minimum number of $k$--linear forests (forests whose connected components are paths of length at most $k$) required to partition the edge set of $G$, and it is denoted by $la_k(G)$.
Bermond et al.~\cite{BerFouHabPer} conjectured that $la_5(G)=2$ for every cubic graph $G$; in other words that the edge set of a cubic graph can be partitioned into two $k$--linear forests for all $k\geq5$. The conjecture is proved by Thomassen in~\cite{Tho}.

In the same paper, Thomassen remarks that $5$ cannot be replaced by $4$ because of the two cubic graphs of order $6$, but up until now it was unknown if there exists a graph with at least eight vertices for which $5$ cannot be replaced by $4$.

In fact, Jackson and Wormald explicitly ask in~\cite{JakWor} if it is true that every cubic graph of order at least eight can be decomposed in two $4$--linear forests.

\begin{problem}[Jackson and Wormald~\cite{JakWor}]
Is it true that $la_4(G)=2$ for all cubic graphs $G$ with at least eight vertices?
\end{problem}

\begin{figure}
\centering
\includegraphics[width=6cm]{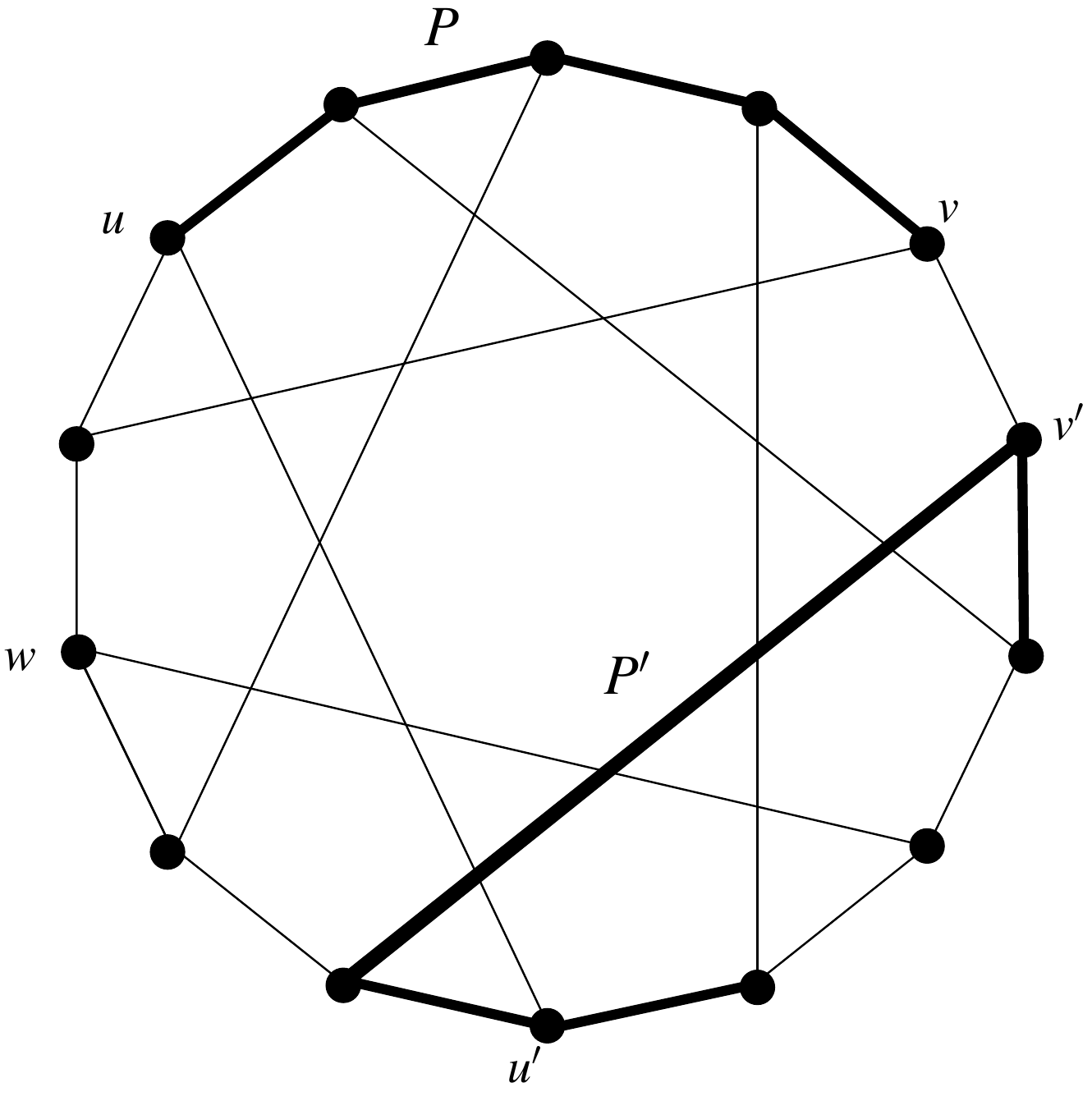}
\caption{The paths $P$ and $P'$ of Heawood graph used in the proof of Theorem~\ref{Th:Heawood}}\label{fig:heawood}

\end{figure}

We give in Theorem~\ref{Th:Heawood} a negative answer to this question by showing that the Heawood graph (see Figure~\ref{fig:heawood}) cannot be decomposed in two $k$--linear forests for $k<5$. We also computationally show that the two cubic graphs of order 6 and the Heawood graph are the only graphs with this property up to at least 28 vertices (cf.\ Observation~\ref{obs:linear_arboricity}).

\begin{lemma}\label{lem:av_three}
Let $G$ be a cubic graph. Consider a decomposition of $G$ into two linear forests. Then, the average length of a path in the decomposition is three.
\end{lemma}

\begin{proof}
An internal vertex of a path is a vertex of degree $2$ in the path. Similarly, an external vertex is a vertex of degree $1$ in the path.
Observe that every vertex of $G$ is internal in one linear forest and external in the other one.
Since the number of paths in the decomposition is $|V(G)|/2$ (every path has exactly two external vertices), their average length is $\frac{2|E(G)}{|V(G)|}=3$.
\end{proof}

Recall that the Heawood graph $H$ is the point/line incidence graph of the Fano plane $PG(2,2)$ which is bipartite, $3$--regular, $3$--arc transitive and has girth~$6$ and diameter~$3$.

\begin{theorem}\label{Th:Heawood}
Let $H$ be the Heawood graph. Then, $la_4(H)>2$.
\end{theorem}

\begin{proof}
First of all, note that $H$ cannot have a decomposition into two $3$--linear forests. Indeed, by Lemma~\ref{lem:av_three}, if the two linear forests have all paths of length at most three, then they have all paths of length exactly three. Moreover, since all paths have the same length, the number of paths in each forest is the same. A contradiction, since we have an even number of paths with three edges, but $H$ has $21$ edges, which is not a multiple of $6$.

Now, suppose that $H$ has a decomposition into two $4$--linear forests, say $F_1$ and $F_2$, and by the previous observation we can assume that at least one of them has a path of length $4$, say $P \in F_1$.
Consider the two external vertices of $P$, say $u$ and $v$. The vertices $u$ and $v$ are both internal vertices in paths of $F_2$. Since $P$ has even length, $u$ and $v$ are in the same bipartition class. Moreover, $H$ has diameter $3$, then $u$ and $v$ are at distance $2$, and the four edges in $F_2$ which are incident to them are all distinct and form a path $Q \in F_2$ of length $4$ . Now, we consider the external vertices of $Q$, say $u'$ and $v'$. Again, the four edges incident $u'$ and $v'$ and not in $Q$ form a further path of length $4$, say $P' \in F_1$. Since $H$ is $3$--arc transitive, the paths $Q$ and $P'$ are, up to symmetries, independent by the choice of the first path $P$ (see Figure~\ref{fig:heawood}).
There is no way to add other paths to $P$ and $P'$ in order to obtain a linear forest since $H-V(P\cup P')$ is a copy of $K_{1,3}$ which does not admit a spanning linear forest (the degree $3$ vertex of $K_{1,3}$ is denoted by $w$ in Figure~\ref{fig:heawood}).
\end{proof}

We also wrote a program to test if a cubic graph can be decomposed in two 4--linear forests. The source code of this program can be downloaded from~\cite{bisection-site}. Using the generator for cubic graphs \textit{snarkhunter}~\cite{snarkhunter-site,brinkmann_11}, we generated all cubic graphs up to 28 vertices and then tested which of these graphs cannot be decomposed in two 4--linear forests. This led to the following observation:

\begin{observation} \label{obs:linear_arboricity}
The two cubic graphs of order $6$ and the Heawood graph are the only cubic graphs up to at least $28$ vertices which cannot be decomposed in two $4$--linear forests.
\end{observation}

Following the flavour of what we have discussed in the previous sections and, particularly, about Wormald's Conjecture, we wonder how we have to modify Thomassen's result if we require to have two linear forests which are also isomorphic. Obviously, this only makes sense if we consider cubic graphs of order $0 \pmod 4$. We have extended our program to take this into account and surprisingly, it turns out that we do not need paths of length more than four for all graphs of order at most 24:

\begin{observation}
Every cubic graph with $|V| \equiv 0 \pmod 4$ with at most 24 vertices has a decomposition into two isomorphic $4$--linear forests.
\end{observation}

We leave the following strengthened version of Wormald's Conjecture as an open problem:

\begin{problem}\label{pro:4linearforests}
Does every cubic graph with $|V| \equiv 0 \pmod 4$ have a decomposition into two isomorphic $4$--linear forests?
\end{problem}

Another interesting and related open question posed by Jackson and Wormald in~\cite{JakWor} is the following:

\begin{problem}[Jackson and Wormald~\cite{JakWor}]
Does every cubic graph have a $2$--edge colouring such that each monochromatic connected component has at most four edges?
\end{problem}

Note that a connected component is not necessarily a path here and a direct check shows that the Heawood graph is not a counterexample in this case.

Furthermore, the two cubic graphs on 6 vertices admit a $2$--edge colouring such that each monochromatic connected component has at most four edges. In particular, the bipartite cubic graph on 6 vertices admits, up to symmetry, only one $2$--edge colouring with all monochromatic connected component having at most four edges, and exactly one of these connected components has precisely four edges.

Hence, it follows that we cannot have three as an upper bound in the previous problem and, from Observation~\ref{obs:linear_arboricity}, that the smallest counterexample must have at least 30 vertices.

\section*{Acknowledgements}
We wish to thank the anonymous referees for their valuable suggestions.
Most computations for this work were carried out using the Stevin Supercomputer Infrastructure at Ghent University.







\section*{Appendix: Cubic Graphs with no Strong Wormald Colouring}\label{sect:strong_wormald}

\begin{figure}[h]
\centering
\includegraphics[width=10cm]{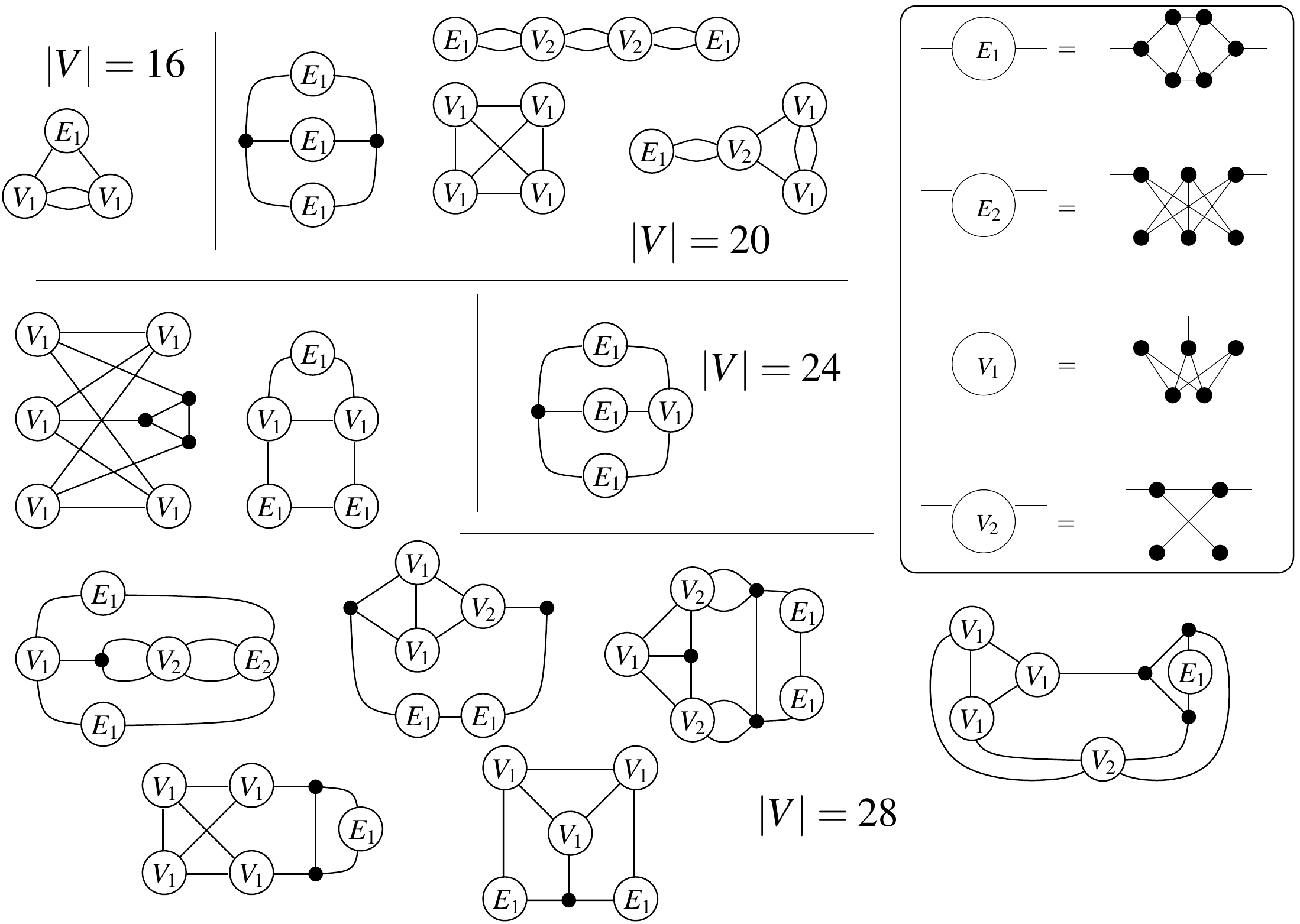}
\caption{A complete list of all $2$--connected cubic graphs without a Strong Wormald Colouring of order at most 28}\label{fig:2connectedwithoutstrongwormald}
\end{figure}

According to Conjecture~\ref{Wconj} every cubic graph of order a multiple of four admits a Wormald Colouring.
We have defined a Strong Wormald Colouring as a Wormald Colouring of the edges of a cubic graph which avoids paths of length one.
We have also observed that cubic graphs with a Strong Wormald Colouring cannot be counterexamples for Ando's Conjecture. Hence,
it is natural to search for cubic graphs (of order a multiple of four) without a Strong Wormald Colouring. We have completed an exhaustive search for cubic graphs without a Strong Wormald Colouring up to 28 vertices (see Section~\ref{comput:wormald} for a description of the algorithm and graph counts) and were able to find such graphs: the smallest one is the graph with 16 vertices at the top of Figure~\ref{fig:2connectedwithoutstrongwormald}.
However, at least for the $2$-connected case, we have found very few examples of cubic graphs without a Strong Wormald Colouring and they share some common features. In particular, note that all $2$-connected examples admit a $3$--edge colouring and hence they admit a $2$--bisection and an Ando Colouring as well.

In Figure~\ref{fig:2connectedwithoutstrongwormald}, we propose a compact description of all $2$--connected cubic graphs without a Strong Wormald Colouring up to 28 vertices. All of them represent several copies of the complete bipartite graph $K_{3,3}$ of order 6 with some edges/vertices removed. In particular, the subgraphs $E_i$ and $V_i$ ($i=1,2$) are obtained by removing $i$ edges and vertices, respectively, from a copy of $K_{3,3}$.

\begin{figure}[h]
\centering
\includegraphics[width=10cm]{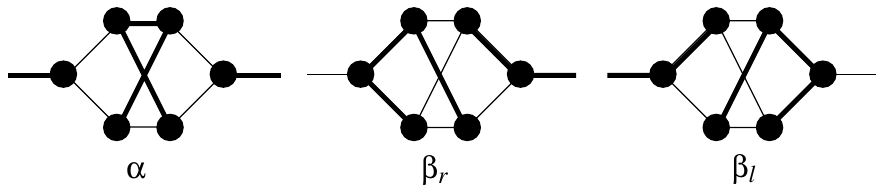}
\caption{Different types of colourings of $E_1$}\label{fig:typesE1}
\end{figure}

To conclude this paper, we show how we can generalise some of the examples from Figure~\ref{fig:2connectedwithoutstrongwormald} in order to obtain an infinite family of graphs without a Strong Wormald Colouring. First of all, we observe that a Strong Wormald Colouring induces, up to symmetry, only two different types of edge colourings of the subgraph $E_1$. We will call these two induced colourings type $\alpha$ and type $\beta$ colourings (see Figure~\ref{fig:typesE1}).
In particular,  we represent in Figure~\ref{fig:typesE1} a colouring of type $\alpha$ and two symmetric colourings of type $\beta$. We denote by $\beta_l$ and $\beta_r$ a colouring of type $\beta$ where the pending edge which belongs to a path having another edge in the copy of $E_1$ is on the left or on the right, respectively.

\begin{figure}[h]
\centering
\includegraphics[width=7cm]{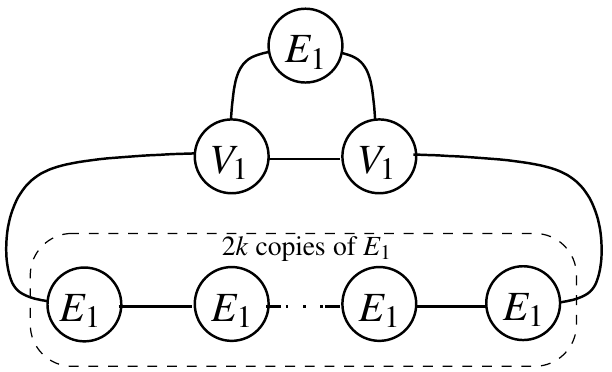}
\caption{The graph $G_k$ }\label{fig:infinite_family}
\end{figure}

Now we construct the infinite family of cubic graphs $G_k$, $k\geq1$ (see Figure~\ref{fig:infinite_family}).
The graph $G_1$ has $28$ vertices and so it appears also in Figure~\ref{fig:2connectedwithoutstrongwormald}.
The graph $G_k$ is obtained starting from $G_1$ by replacing the two copies of $E_1$ with $2k$ consecutive copies of $E_1$.

\begin{proposition}
For every $k\geq 1$, the graph $G_k$ does not admit a Strong Wormald Colouring.
\end{proposition}
\begin{proof}
We argue by induction on $k$. We need to consider the second member of the family, $G_2$, as base case. The graph $G_2$ has $40$ vertices and we have directly checked it does not admit a Strong Wormald Colouring using a computer algorithm (see Section~\ref{comput:wormald} for details of the algorithm).
Now, we assume that $G_{k}$, $k\geq 2$, has no Strong Wormald Colouring and we prove it for $G_{k+1}$.
By contradiction, suppose that $G_{k+1}$ admits a Strong Wormald Colouring. Since $k\geq2$, then $2(k+1)\geq 6$, that is: we have at least $6$ consecutive copies of $E_1$ in $G_{k+1}$.

Assume that one of them, say $E$, is of type $\alpha$. Since the two pending edges of $E$ cannot be in a path of length one, the copy on the left of $E$ is  of type $\beta_r$ and the one on the right of type $\beta_l$. More in general, repeating the same argument, all copies on the left of $E$ are of type $\beta_r$ and all copies on the right of $E$ are of type $\beta_l$. Since we have at least six copies of $E_1$, either we have three consecutive copies of $\beta_r$ or three consecutive copies of $\beta_l$. Either way, we can remove the last two of these three consecutive copies in order to obtain a graph isomorphic to $G_k$. Moreover, there is exactly one path of length two, one path of length three and one path of length five in each monochromatic component of the two removed copies of $E_1$. The two monochromatic subgraphs are still isomorphic in $G_k$ and all paths are still of length at least two, hence the graph $G_k$ inherits a Strong Wormald Colouring, a contradiction.

Therefore, assume no copy of $E_1$ is of type $\alpha$. On the right of a copy of type $\beta_r$ there must be another copy of type $\beta_r$, and the same holds for copies of type $\beta_l$. Hence, since we have at least six copies of $E_1$, we have three consecutive copies of type $\beta_r$ or three consecutive copies of type $\beta_l$. Either way we obtain a contradiction by the same argument used above.
\end{proof}

Finally, note that $G_k$ is  a $3$-edge-colourable (bipartite) graph for every $k$. Hence, $G_k$ admits a $2$--bisection and by Proposition~\ref{BanLinialimpliesAndo} is not a counterexample for Ando's Conjecture.

\end{document}